\def\wlog#1{}
 \makeatletter
 \def\@latex@info#1{}
 \def\@font@info#1{}
 \makeatother

\documentclass[mathsec, numreferences, 12pt]{kluwer}

\usepackage{upref}
\usepackage{amsmath}

\newcommand{\oldendpf}{}
\let\oldendpf=\endpf
\renewcommand{\endpf}{\qed\oldendpf}

\newcommand{\equationlabel}{\label}

\makeatletter

\@openrightfalse

\renewcommand{\labelitemi}{\m@th$\bullet$}
\renewcommand{\labelitemii}{\m@th$-$}

\@noidtrue

\renewcommand{\equationlabel}[1]{%
 \def\@currentlabel{\theequation}
 \protect\ilabel{\thearticle #1}}

\newcommand{\proofname}{Proof}
\renewcommand{\pf}{\par 
    \addvspace{1\baselineskip 
            \@plus 0.5\baselineskip \@minus 0.1\baselineskip}%
    \indent
    {\it \proofname.\/} \ignorespaces }

\def\kap@enumerate[#1]{%
     \ifnum \@enumdepth >3 \@toodeep\else
     \advance\@enumdepth \@ne
     \edef\@enumctr{enum\romannumeral\the\@enumdepth}
     \list{\csname label\@enumctr\endcsname}{%
       \kapenumargs 
       \usecounter{\@enumctr}
       \settowidth\labelwidth{#1.}
       \setlength{\leftmargin}{\labelwidth} 
       \addtolength{\leftmargin}{\labelsep}
       \addtolength{\leftmargin}{\parindent} 
       \def\makelabel##1{\hss \llap{##1}}}%
     \fi
   }

\def\kap@itemize[#1]{\def\klu@arg{#1}%
    \ifnum \@itemdepth >3 \@toodeep
    \else
      \advance\@itemdepth \@ne 
      \edef\@itemitem{labelitem\romannumeral\the\@itemdepth}%
      \ifx \klu@arg\empty 
        \list {\csname\@itemitem\endcsname}%
        {\kapitemargs
         \def\makelabel##1{\hss \llap{##1}}}%
      \else 
        \list {\klu@arg }%
        {\kapitemargs
         \def\makelabel##1{\hss \llap{##1}}}%
      \fi
    \fi 
    }

\def\@opargbegintheorem#1#2#3{\trivlist 
      \global\@novspacetrue \itemindent\theoremsep 
      \item[\kern \labelsep {\@stylehead\@thmscase{#1}\ #2\ 
      \rm(#3). \/}]\ \@styletext}

\def\@opargbegindisplay#1#2#3{\trivlist
      \global\@novspacefalse
      \itemindent \dispsep
      \item[{\@disphead \@dispcase{#1}\ #2\/\ 
      \upshape({#3}).}]\@disptext}

\def\@stylehead{\bf }
\def\@thmscase{}

\makeatother

\newtheorem{thm}[equation]{Theorem}
\newtheorem{lem}[equation]{Lemma}
\newtheorem{cor}[equation]{Corollary}

\newdisplay{defn}[equation]{Definition}
\newdisplay{eg}[equation]{Example}
\newdisplay{notation}[equation]{Notation}
\newdisplay{rem}[equation]{Remark}

\newcommand{\oldsection}{}
 \let\oldsection=\section
\renewcommand{\section}{\setcounter{equation}{0}\oldsection}
\renewcommand{\theequation}{\thesection.\arabic{equation}}

\newenvironment{thmref}{\thmrefer}{}
\newcommand{\thmrefer}[1]{\renewcommand\theequation
  {\protect\ref{#1}$'$}\addtocounter{equation}{-1}}

 \newcounter{case}
 \newenvironment{case}[1][\unskip]{\refstepcounter{case}
 \setcounter{subcase}{0}\em
 \medskip \noindent Case \thecase\ #1.\ }{\unskip\upshape}
 \renewcommand{\thecase}{\arabic{case}}

 \newcounter{subcase}
 \newenvironment{subcase}[1][\unskip]{\refstepcounter{subcase}\em
 \medskip \noindent Subcase \thesubcase\ #1.\ }{\unskip\upshape}
 \renewcommand{\thesubcase}{\thecase.\arabic{subcase}}

 \newenvironment{Claim}[1][\unskip]{\em
 \medskip \noindent Claim\ #1.\ }{\unskip\upshape}

 \newcounter{step}

\newenvironment{step}[1][\unskip]{\refstepcounter{step}\em
 \medskip \noindent Step \thestep\ #1.\ }{\unskip\upshape}
 \renewcommand{\thestep}{\arabic{step}}

 \newcounter{stepclaim}

\newenvironment{stepclaim}[1][\unskip]{\refstepcounter{stepclaim}\em
 \medskip \noindent Claim \thestepclaim\ #1.\
}{\unskip\upshape}
 \renewcommand{\thestepclaim}{\thestep.\arabic{stepclaim}}

\newcommand{\pref}[1]{{\upshape(}\ref{#1}{\upshape)}}
\newcommand{\see}[1]{{\upshape(}see~\ref{#1}{\upshape)}}
\newcommand{\cf}[1]{{\upshape(}cf.~\ref{#1}{\upshape)}}
\newcommand{\fullref}[2]{\ref{#1}\pref{#1-#2}}

\newcommand{\Lie}[1]{\mathfrak{\lowercase{#1}}}
\newcommand{\sgn}{\operatorname{sgn}}
\newcommand{\supp}{\operatorname{supp}}
\newcommand{\ad}{\operatorname{ad}}
\newcommand{\Ad}{\operatorname{Ad}}
\newcommand{\Rad}{\operatorname{Rad}}
\newcommand{\SL}{\operatorname{SL}}
\newcommand{\GL}{\operatorname{GL}}
\newcommand{\PSL}{\operatorname{PSL}}
\newcommand{\SO}{\operatorname{SO}}
\renewcommand{\Sp}{\operatorname{Sp}}
\newcommand{\SU}{\operatorname{SU}}
\newcommand{\Id}{\operatorname{Id}}
\newcommand{\real}{\mathord{\mathbb{R}}}
\newcommand{\complex}{\mathord{\mathbb{C}}}
\newcommand{\rational}{\mathord{\mathbb{Q}}}

\newcommand{\iso}{\cong}
\newcommand{\funddom}{\sigma}
\newcommand{\Homeo}{\operatorname{Homeo}}
\newcommand{\HomeoLeb}
  {\operatorname{Homeo}^{\text{\upshape Leb}}}
\newcommand{\Diff}{\operatorname{Diff}}
\newcommand{\Prob}{\operatorname{Prob}}
\newcommand{\Frank}[1]{\mathop{\mbox{$#1${\upshape-rank}}}}
\newcommand{\Rrank}{\Frank{\real}}
\newcommand{\Rot}{\operatorname{Rot}}
\newcommand{\Isom}{\operatorname{Isom}}
\newcommand{\algG}{\mathord{\mathbf{G}}}

\newcommand{\algA}{\mathord{\mathbf{A}}}
\newcommand{\algU}{\mathord{\mathbf{U}}}
\newcommand{\torus}{\mathbb{T}}
\newcommand{\atimes}{\mathbin{\times_\alpha}}
\newcommand{\ttimes}{\mathbin{\times_\tau}}
\newcommand{\btimes}{\mathbin{\times_\beta}}
\newcommand{\Sl}{\operatorname{\Lie{sl}}}
\newcommand{\invariant}{\mathord{\mathcal{F}}}

\newcommand{\bigset}[2]{\left\{\, #1 
 \mathrel{\left| \vphantom {\left\{ #1 \mid #2 \right\} } \right.}
 #2 \,\right\} }

\begin{document}
\begin{article}
\begin{opening}

\title{Actions of semisimple Lie groups on circle bundles}

\author{Dave \surname{Witte}
 \email{dwitte@math.okstate.edu,
  http://www.math.okstate.edu/\char'176dwitte}} 
 \institute{Department of Mathematics, Oklahoma State
University, Stillwater, OK 74078, USA}

\author{Robert J.~\surname{Zimmer}
 \email{r-zimmer@uchicago.edu}}
 \institute{Department of Mathematics, University of Chicago,
Chicago, IL 60637, USA}

\runningtitle{Actions on circle bundles}

\begin{abstract} 
 Suppose $G$ is a connected, simple, real Lie group with
$\Rrank(G) \ge 2$, $M$~is an ergodic $G$-space with
invariant probability measure~$\mu$, and $\alpha \colon G
\times M \to \Homeo(\torus)$ is a Borel cocycle.
 We use an argument of \'E.~Ghys to show that there is a
$G$-invariant probability measure~$\nu$ on the skew product
$M \atimes \torus$, such that the projection of~$\nu$ to~$M$
is~$\mu$. Furthermore, if $\alpha(G \times M) \subset
\Diff^1(\torus)$, then $\nu$ can be taken to be equivalent
to $\mu \times \lambda$, where $\lambda$ is Lebesgue
measure on~$\torus$; therefore, $\alpha$ is cohomologous to
a cocycle with values in the isometry group of~$\torus$.
 \end{abstract}

\date{\today} 

\end{opening}

\section{Introduction} \label{intro}

\'E.~Ghys \cite{Ghys} recently proved that irreducible
lattices in most semisimple Lie groups of higher real rank
do not have any interesting differentiable actions on the
circle~$\torus$.

\begin{defn} \label{irredlattice}
 A lattice~$\Gamma$ in a connected, semisimple, real Lie
group~$G$ is \emph{irreducible} if $N \Gamma$ is dense
in~$G$, for every closed, connected, noncompact, normal
subgroup~$N$ of~$G$.
 \end{defn}

\begin{notation}
 We use $\Diff^1(\torus)$ to denote the group of
$C^1$~diffeomorphisms of~$\torus$, and $\Diff^1_+(\torus)$
to denote the subgroup of orientation-preserving
diffeomorphisms.
 \end{notation}

\begin{thm}[{Ghys \cite[Thm.~1.1]{Ghys}}]
\label{GhysC1LatticeThm}
 Let $\Gamma$ be an irreducible lattice in a connected,
semisimple, real Lie group~$G$, such that
 \begin{enumerate}
 \item $\Rrank G \ge 2$; and
 \item \label{GhysC1LatticeThm-noSL2}
 there is no continuous homomorphism from~$G$ onto
$\PSL(2,\real)$.
 \end{enumerate}
 Then every homomorphism from $\Gamma$ to~$\Diff^1(\torus)$
has finite image.
 \end{thm}

\begin{rem} \label{GhysLatticeGeneralizations}
 Under the additional assumption that $H^2(\Gamma;\real) =
0$ (and in many other cases), the conclusion of the
theorem was also proved by M.~Burger and N.~Monod
\cite{BurgerMonod-1, BurgerMonod-2, BurgerMonod-3}, as a
consequence of vanishing theorems for bounded cohomology.
(The results of Burger and Monod also apply to the setting
where $\real$ is replaced by other local fields; for
example, $\Gamma$ could be an $S$-arithmetic group
(cf.~\ref{Ghys-Slattice} and~\ref{Sarith-on-circle}).) For a more
restricted class of lattices in real semisimple Lie groups,
B.~Farb and P.~Shalen \cite{FarbShalen} proved finiteness of the
image of homomorphisms into the group $\Diff^{\omega}(M)$ of real
analytic diffeomorphisms of some higher-dimensional
manifolds.
 \end{rem}

In this paper, we extend Ghys' Theorem to the context of
semisimple Lie group actions on circle bundles, or, more
generally, $\Diff^1(\torus)$-valued Borel cocycles for
ergodic actions of~$G$. We first recall:

\begin{defn}[{\cite[Defns.~4.2.1 and~4.2.2, p.~65, and top
of p.~75]{ZimmerBook}}]
 Suppose $M$ is a Borel $G$-space with quasi-invariant
measure~$\mu$, and $H$ is a topological group (such that
the Borel structure on~$H$ is countably generated). 
 \begin{itemize}
 \item A Borel function $\alpha \colon G \times M \to H$ is
a \emph{Borel cocycle} if, for all $g,h \in G$, we have
$\alpha(gh, m) =\alpha(g, hm) \, \alpha(h,m)$ for a.e.~$m
\in M$.
 \item Two Borel cocycles $\alpha, \beta \colon G \times M
\to H$ are \emph{cohomologous} if there is a Borel function
$\phi \colon M \to H$, such that, for each $g \in G$, we
have $\beta(g,m) = \phi(gm)^{-1} \alpha(g,m) \phi(m)$, for
a.e.~$m \in M$.
 \item A Borel cocycle $\alpha \colon G \times M \to H$ is
\emph{strict} if, for all $g,h \in G$, we have
$\alpha(gh, m) =\alpha(g, hm) \, \alpha(h,m)$ for
\emph{every} $m \in M$. For every Borel cocycle $\alpha
\colon G \times M \to H$, there is a strict Borel cocycle
$\alpha' \colon G \times M \to H$, such that, for every $g
\in G$, we have $\alpha'(g,m) = \alpha(g,m)$ for a.e.\ $m
\in M$ \cite[Thm.~B.9, p.~200]{ZimmerBook}.
 \item If $\alpha \colon G \times M \to H$ is a strict Borel
cocycle and $S$ is a Borel $H$-space, the
\emph{skew-product} action $M \atimes S$ is the Borel
action of~$G$ on $M \times S$ defined by
 $g\cdot (m,s) = \bigl( gm, \alpha(g,m) s \bigr)$.
 \end{itemize}
 \end{defn}

 Recall that any smooth action on a circle bundle defines a
$\Diff^1(\torus)$-valued cocycle on the base, and that the
action on the bundle is measurably conjugate to the skew
product action defined by this cocycle.  Conversely, the
skew product defined by any $\Diff^1(\torus)$-valued cocycle
can be viewed as an action on a measurable circle bundle
over the base.

For $M = G/\Gamma$, cohomology classes of Borel cocycles
$\alpha \colon G \times M \to \Diff^1(\torus)$ are in
bijective correspondence with conjugacy classes of
homomorphisms $\hat\alpha \colon \Gamma \to \Diff^1(\torus)$
\cite[Prop.~4.2.13, p.~70]{ZimmerBook}. Then the conclusion
of Ghys' Theorem asserts that $\alpha$ is cohomologous to a
Borel cocycle whose image is a finite subgroup of
$\Diff^1(\torus)$. However, the following example shows
that this conclusion is not valid for Borel cocycles for
more general $G$-spaces; not even for Borel cocycles that
arise from a $C^\infty$, volume-preserving action of~$G$ on
a principal $\torus$-bundle over a compact manifold.

\begin{eg}
 Let 
 \begin{itemize}
 \item $H$ be a connected, semisimple Lie group;
 \item $\Gamma$ be a torsion-free, cocompact
lattice in~$H$;
 \item $T$ be a subgroup of~$H$ that is isomorphic
to~$\torus$;
 \item $G$ be a closed subgroup of~$H$ that
centralizes~$T$ and acts ergodically on $H/\Gamma$
\see{MooreErgodicityThm}; and 
 \item $M = T \backslash H/\Gamma$.
 \end{itemize}
 Because $\Gamma$ is torsion free and cocompact, we know
that $M$ is a compact manifold. Because $G$
centralizes~$T$, the action of~$G$ by translation on
$H/\Gamma$ factors through to an action on~$M$; we see that
$H/\Gamma$ is a principal $\torus$-bundle over~$M$, and $G$
acts on~$H/\Gamma$ by bundle automorphisms. Thus, there is
a Borel cocycle $\alpha \colon G \times M \to \torus$, such
that the action of~$G$ on~$H/\Gamma$ is isomorphic to the
skew product $M \atimes \torus$. By assumption, the action
of $G$ on $H/\Gamma$ is ergodic, so, if $\beta$ is any
cocycle cohomologous to~$\alpha$, then $M \btimes \torus$
must be ergodic. Therefore, the image of~$\beta$ cannot be
contained in any finite group of transformations
of~$\torus$.
 \end{eg}

These examples show that there can be nontrivial cocycles
into $\Isom(\torus)$, the isometry group of~$\torus$. Our
extension of Ghys' Theorem shows that if $G$ has Kazhdan's
property~$(T)$ \see{KazhdanDefn}, then every cocycle into
$\Diff^1(\torus)$ for a much more general $G$-action is
cohomologous to one into $\Isom(\torus)$. (However, as far
as we know, the homeomorphisms in the image of the map
implementing the cohomology may not be differentiable, but
only absolutely continuous.) In more geometric terms, this
asserts that for $G$-actions on very general circle
bundles, there is a measurable choice of metric on each
fiber that is preserved by the action. I.e., the action on
the bundle is an ``isometric extension" of the base.

\begin{defn} \label{irredaction}
 Let $G$ be a connected, semisimple Lie group, and let $M$
be an ergodic $G$-space with quasi-invariant measure. We
say that $M$ is \emph{irreducible} if every closed,
connected, noncompact, normal subgroup of~$G$ is ergodic
on~$M$.
 \end{defn}

\begin{notation}
 $\HomeoLeb(\torus)$ denotes the group of all
homeomorphisms~$\phi$ of~$\torus$, such that
$\phi_*\lambda$ has the same null sets
as~$\lambda$, where $\lambda$ is the Lebesgue measure
on~$\torus$.
 \end{notation}

\begin{thmref}{C1-isometry}
 \begin{thm}
 Let 
 \begin{itemize}
 \item $G$ be a connected, real, semisimple Lie group, such that
 \begin{itemize}
 \item $G$ has Kazhdan's property~$(T)$, and 
 \item $\Rrank G \ge 2$;
 \end{itemize}
 \item $M$ be an irreducible ergodic $G$-space with finite
invariant measure~$\mu$; and
 \item $\alpha \colon G \times M \to \Diff^1(\torus)$ be a Borel
cocycle.
 \end{itemize}
 Then, as a cocycle into $\HomeoLeb(\torus)$, $\alpha$
is cohomologous to a cocycle with values in $\Isom(\torus)$.

 Furthermore, if $\alpha(g,m)$ is orientation preserving,
for almost every $(g,m) \in G \times M$, then, as a cocycle
into $\HomeoLeb(\torus)$, $\alpha$~is cohomologous to a
cocycle with values in the rotation group $\Rot(\torus)$.
 \end{thm}
 \end{thmref}

It is an open question whether Ghys'
Theorem~\ref{GhysC1LatticeThm} remains valid if
$\Diff^1(\torus)$ is replaced with the homeomorphism group
$\Homeo(\torus)$. (Witte~\cite{Witte-circle} showed that
the answer is affirmative if $\Gamma$ is an arithmetic
lattice of $\rational$-rank at least two.) However, Ghys
(and, in most cases, also Burger and Monod) made the
following major step toward an affirmative answer.

\begin{thm}[{Ghys, cf.\ \cite[Thm.~3.1]{Ghys}}]
\label{GhysHomeoLattThm}
 Let $\Gamma$ be an irreducible lattice in a connected,
semisimple, real Lie group~$G$, such that
 \begin{enumerate}
 \item $\Rrank G \ge 2$; and
 \item \label{GhysThm-noSL2}
 there is no continuous homomorphism from~$G$ onto
$\PSL(2,\real)$.
 \end{enumerate}
 Then every continuous action of~$\Gamma$ on~$\torus$ has
an invariant probability measure.

 In fact, every continuous action of~$\Gamma$ on~$\torus$
has a finite orbit.
 \end{thm}

Ghys obtained Theorem~\ref{GhysC1LatticeThm} by combining
Theorem~\ref{GhysHomeoLattThm} with the Thurston Stability
Theorem~\ref{ReebThurston}. (He also proved that if $G$
does have a continuous homomorphism onto $\PSL(2,\real)$,
then any action of~$\Gamma$ on~$\torus$ either preserves a
probability measure or is semi-conjugate to a finite cover
of the restriction of a $G$-action \cf{Ghys-SL2}.)

\begin{thmref}{ReebThurston}
 \begin{thm}[{Thurston \cite{Thurston}}]
 Suppose $\Gamma$ is a finitely generated group, such that
$\Gamma/[\Gamma,\Gamma]$ is finite. If $\sigma \colon
\Gamma \to \Diff^1_+(\torus)$ is any homomorphism, such
that $\sigma(\Gamma)$ has a fixed point, then
$\sigma(\Gamma)$ is trivial.
 \end{thm}
 \end{thmref}

The following theorem is the natural generalization of
Theorem~\ref{GhysHomeoLattThm} to the setting of ergodic
$G$-actions. Although Ghys did not state this result, it can
be proved by translating his proof in a straightforward way
from the setting of homomorphisms of lattices to the
setting of Borel cocycles for ergodic $G$-actions. In
Section~\ref{ProveGhys}, we provide a proof that is based
on Ghys' ideas, but is much shorter than a direct
translation.

\begin{thm} \label{GhysThm}
 Let 
 \begin{itemize}
 \item $G$ be a connected, semisimple, real Lie group, such that
 \begin{itemize}
 \item $\Rrank G \ge 2$, and
 \item there is no continuous homomorphism from~$G$ onto
$\PSL(2,\real)$;
 \end{itemize}
 \item $M$ be an irreducible ergodic $G$-space with
invariant probability measure~$\mu$; and
 \item $\alpha \colon G \times M \to \Homeo(\torus)$ be a
strict Borel cocycle.
 \end{itemize}
 Then there is a $G$-invariant probability measure~$\nu$ on
$M \atimes \torus$, such that the projection of~$\nu$
to~$M$ is~$\mu$.
 \end{thm}

We obtain Theorem~\ref{C1-isometry} by combining
Theorem~\ref{GhysThm} with the following generalization of
Theorem~\ref{ReebThurston}.

\begin{defn}
 Let $\alpha \colon G \times M \to H$ be a Borel cocycle,
and let $Y$ be an $H$-space. A function $f \colon M \to Y$
is \emph{$\alpha$-equivariant} if, for each $g \in G$, we
have $f(gm) = \alpha(g,m) f(m)$ for almost every $m \in M$.
 \end{defn}

 \begin{thmref}{C1-fp-measure}
 \begin{thm}
 Let 
 \begin{itemize}
 \item $G$ be a connected Lie group with Kazhdan's
property~$(T)$;
 \item $M$ be an ergodic $G$-space with finite invariant
measure~$\mu$; 
 \item $\alpha \colon G \times M \to \Diff^1(\torus)$ be a Borel
cocycle; and
 \item $f \colon M \to \torus$ be an
$\alpha$-equivariant measurable map. {\upshape(}In
bundle theoretic terms, $f$ is a measurable
$G$-invariant section.{\upshape)}
 \end{itemize}
 Then, as a cocycle into $\HomeoLeb(\torus)$,
$\alpha$~is cohomologous to the trivial cocycle.
 \end{thm}
 \end{thmref}

Theorems~\ref{C1-isometry} and~\ref{GhysThm} can be
generalized to allow~$G$ to be a $S$-algebraic group
\see{Ghys-S-algebraic}, and there are also analogues for
$\Gamma$-actions, where $\Gamma$ is a lattice in~$G$
\see{Ghys-Lattice}. Thus, as was already mentioned in
Remark~\ref{GhysLatticeGeneralizations}, Ghys'
Theorem~\ref{GhysC1LatticeThm} can be generalized to allow
$\Gamma$ to be an $S$-arithmetic group
\see{Sarith-on-circle}.

The paper is organized as follows.
 Section~\ref{prelims} establishes notation, and recalls
various results from measure theory, Lie theory, ergodic
theory, and Kazhdan's property~$(T)$. 
 Section~\ref{crucial-lemma} constructs a pair of subgroups
that play a crucial role in the proof of
Theorem~\ref{GhysThm}, which is presented in
Section~\ref{ProveGhys}.
 Section~\ref{Thurston-section} proves
Theorems~\ref{C1-fp-measure} and~\ref{C1-isometry}, our
results on differentiable actions.
 Section~\ref{versions-section} extends our main results to
slightly different settings.

\begin{acknowledgements}
 This research was partially supported by grants from the
National Science Foundation (DMS--9801136 and DMS--9705712).
 D.W.\ is grateful to \'Etienne Ghys for many very helpful
discussions about lattice actions on the circle, and to
M.~Burger and Y.~Shalom for instructive comments that
greatly improved the results on $S$-arithmetic groups. He
would like to thank the \'Ecole Normale Sup\'erieure de
Lyon, the University of Chicago, and the University of
Bielefeld for their hospitality while this research was
underway, and he would also like to thank the
German-Israeli Foundation for Research and Development for
financial support that made the visit to Bielefeld
possible.
 \end{acknowledgements}

\section{Preliminaries} \label{prelims}

\subsection{Probability measures}

\begin{notation}
 We use $I$ to denote the unit interval $[0,1]$, and
$\torus$ to denote the unit circle.
 For $\Omega
= \torus$ or~$I$:
 \begin{itemize}
 \item $\lambda$ denotes the Lebesgue measure on~$\Omega$;
and
 \item $\Prob(\Omega)$ denotes the space of probability
measures on~$\Omega$, with the weak* topology.
 \end{itemize}
 \end{notation}

\begin{defn}
 Measures $\mu_1$ and~$\mu_2$ on a Borel space~$X$ are
\emph{equivalent} (or in the same \emph{measure class}) if
they have the same null sets.
 \end{defn}

\begin{lem}
 Let $\alpha \colon G \times M \to \HomeoLeb(\torus)$ be a
Borel cocycle. There is a $G$-invariant probability measure
on $M \atimes \torus$ that is equivalent to $\mu \times
\lambda$ if and only if $\alpha$ is equivalent to a cocycle
with values in $\Isom(\torus)$.
 \end{lem}

\begin{pf}
 ($\Leftarrow$) By assumption, there is a Borel cocycle
$\beta \colon G \times M \to \Isom(\torus)$, and a Borel
function $\phi \colon M \to \HomeoLeb(\torus)$, such that,
for each $g \in G$, we have 
 $$\alpha(g,m) = \phi(gm)^{-1} \beta(g,m) \phi(m)$$
 for a.e.~$m \in M$. Let
 $$ \nu = \int_M \bigl( m \times \phi(m)^{-1}_* \lambda
\bigr) \, d \mu(m)  \in \Prob(M \times \torus) .$$
 Because $\phi(m) \in \HomeoLeb(\torus)$, we know that
$\phi(m)_* \lambda$ is equivalent to~$\lambda$, for every
$m \in M$, so $\nu$ is equivalent to $\mu \times \lambda$.
Because $\lambda$ is invariant under $\Isom(\torus)$, it is
easy to see that $\nu$ is invariant under the action of~$G$
on $M \atimes \torus$.

($\Rightarrow$) Let $\nu$ be a $G$-invariant probability measure on $M
\atimes \torus$ that is equivalent to $\mu \times \lambda$.
We may write
 $$ \nu = \int_M ( m \times \nu_m ) \, d \mu(m) ,$$
 where $\nu_m$ is a probability measure on~$\torus$.
Because $\nu$ is equivalent to $\mu \times \lambda$, we
know that $\nu_m$ is equivalent to~$\lambda$, for a.e.~$m
\in M$. Thus, for a.e.~$m \in M$, there exists $\phi(m) \in
\HomeoLeb(\torus)$, such that $\nu_m = \phi(m)^{-1}_*
\lambda$. Now define $\beta \colon G \times M \to
\HomeoLeb(\torus)$ by
 $\beta(g,m) = \phi(gm) \alpha(g,m) \phi(m)^{-1}$.
 Then $\mu \times \lambda$ is a $G$-invariant measure on $M
\btimes \torus$, so we see, for each $g \in G$, that
$\beta(g,m)$ preserves~$\lambda$, and hence is in
$\Isom(\torus)$, for a.e.\ $m \in M$.
 \end{pf}

\subsection{Lie theory {\protect\cite[Chap.~1]{WarnerBook}}} 
 Let $G$ be a connected, semisimple, real Lie group.

\begin{notation}
 We use lower-case gothic letters $\Lie G, \Lie H, \Lie P,
\Lie Q$, etc.\ for the Lie algebras of Lie groups $G, H, P,
Q$, etc.
 \end{notation}

\begin{defn}
 A subalgebra~$\Lie A$ of~$\Lie G$ is a \emph{maximal split
toral subalgebra} of~$\Lie G$ if 
 \begin{enumerate}
 \item \label{SplitTorus-abel} $\Lie A$ is abelian;
 \item \label{SplitTorus-diag} $\ad_{\Lie G} a$ is
diagonalizable over~$\real$, for every $a \in \Lie A$; and
 \item $\Lie A$ is maximal, with respect to
\pref{SplitTorus-abel} and~\pref{SplitTorus-diag}.
 \end{enumerate}

A \emph{maximal split torus} of~$G$ is a closed, connected
subgroup~$A$ of~$G$, such that the Lie algebra~$\Lie A$
of~$A$ is a maximal split toral subalgebra of~$G$.
 \end{defn}

\begin{defn}
 Let $A$ be a maximal split torus of~$G$.
 \begin{itemize}
 \item For each linear functional $\alpha \colon \Lie A \to
\real$, we let
 $$ \Lie G_\alpha = \{ \, v \in \Lie G \mid 
 \text{$(\ad_{\Lie G} a)(v) = \alpha(a) v$ for all $a \in
\Lie A$} \, \} .$$
 \item A linear functional $\alpha \colon \Lie A \to \real$
is a \emph{real root} of~$\Lie G$ if $\Lie G_\alpha \neq 0$.
 \item The \emph{relative Weyl group} of~$G$ is
$N_G(A)/C_G(A)$.
 \end{itemize}
 \end{defn}

\begin{defn}
 A subalgebra~$\Lie P$ of~$\Lie G$ is \emph{parabolic} if
$\Lie P \otimes \complex$ contains a maximal solvable
subalgebra of~$\Lie G \otimes \complex$. 

 A subgroup~$P$ of~$G$ is \emph{parabolic} if 
 \begin{itemize}
 \item $\Lie P$ is parabolic and
 \item $P = N_G(\Lie P)$.
 \end{itemize}
 \end{defn}

\begin{rem}[{\cite[Thm.~1.2.4.8, p.~75]{WarnerBook}}]
\label{P-structure}
 If $P$ is any parabolic subgroup of~$G$, then $P$ contains
a maximal split torus~$A$ of~$G$. We have $C_G(A) \subset
P$ and, for any real root~$\alpha$ of~$\Lie G$, we have
either $\Lie G_\alpha \subset \Lie P$ or $\Lie G_{-\alpha}
\subset \Lie P$.
 \end{rem}

\begin{rem} \label{parab-sl2}
 A proper subgroup~$P$ of $\SL(2,\real)$ is parabolic if
and only if $P$ is conjugate to
 $ \begin{pmatrix} * & * \\ 0 & * \end{pmatrix} $.
 \end{rem}

\begin{lem} \label{SL2-almparab}
 Let $P$ be a parabolic subgroup of~$G$, and let $L$ be a
closed, connected subgroup of~$G$ that is locally
isomorphic to $\SL(2,\real)$. If $\Lie P \cap\Lie  L$ is  a
parabolic subalgebra of~$\Lie L$, then $P \cap L$ is a
parabolic subgroup of~$L$.
 \end{lem}

\begin{pf}
 Because $\Lie P \cap\Lie  L$ is  a parabolic subalgebra
of~$\Lie L$, there is a parabolic subgroup~$Q$ of~$L$, such
that $Q^\circ = (P \cap L)^\circ$. We wish to show that $Q
\subset P$.
 By definition, $P$ is the normalizer of~$\Lie P$, so it
suffices to show that every subalgebra of~$\Lie G$
normalized by~$Q^\circ$ is also normalized by~$Q$.

Because $\SL(2,\real)$ is simply connected as an algebraic
group, the adjoint representation of~$L$ on~$\Lie G$ must
factor through either $\SL(2,\real)$ or $\PSL(2,\real)$.
Then, because parabolic subgroups of $\SL(2,\real)$ are
Zariski connected, we conclude that every subalgebra
of~$\Lie G$ normalized by~$Q^\circ$ is also normalized
by~$Q$, as desired.
 \end{pf}

\subsection{Ergodic actions}

\begin{thm}[{``Moore Ergodicity Theorem," cf.\
\cite[Thm.~2.2.15, p.~21]{ZimmerBook}}]
\label{MooreErgodicityThm}
 Let 
 \begin{itemize}
 \item $G$ be a connected, semisimple, real Lie group;
 \item $M$ be an irreducible, ergodic $G$-space with finite
invariant measure; and 
 \item $H$ be a closed subgroup of~$G$, such that $\Ad_G H$
is not precompact.
 \end{itemize}
 Then
 \begin{enumerate}
 \item the action of~$H$ on~$M$ is ergodic; and
 \item the diagonal action of~$G$ on $(G/H) \times M$ is
ergodic.
 \end{enumerate}
 \end{thm}

\begin{cor} \label{MooreCor(G/P)^2}
 Let $M$ be an irreducible, ergodic $G$-space with finite
invariant measure, and let $P$ be a minimal parabolic
subgroup of~$G$. Then the diagonal action of~$G$ on $(G/P)
\times (G/P) \times M$ is ergodic.
 \end{cor}

\begin{defn}
 An action of~$G$ on a space~$X$ is \emph{triply
transitive} if $G$ is transitive on the set of ordered
triples of distinct points of~$X$.
 \end{defn}

We note that if $G$ acts triply transitively on~$X$, then
$X$ has no nontrivial, proper $G$-equivariant quotients.
(In particular, every $G$-equivariant quotient of~$X$ is
triply transitive.) Namely, if there are two distinct
points in the same fiber of a quotient map, then, by double
transitivity, $G$ can move them to two points in different
fibers.  This is impossible if the quotient map is
$G$-equivariant.

\subsection{Kazhdan's property~$(T)$}

\begin{defn}[{Kazhdan, cf.\
\cite[Prop.~III.2.8(A), p.~116]{MargulisBook}}]
\label{KazhdanDefn}
 \ 
 A locally compact group~$G$ has \emph{Kazhdan's
property~$(T)$} if, for every unitary representation~$\rho$
of~$G$ on a Hilbert space~$V$, there is a compact
subset~$C$ of~$G$, such that, for every $\epsilon > 0$ there
exists $\delta = \delta(\epsilon) > 0$, such that
 if $v \in V$ is any vector with the property that
 $$ \| \rho(g) v - v \| \le \delta \|v\|
 \text{ \ for every $g \in C$}, $$
 then there is a $\rho(G)$-invariant vector $w \in V$, such
that $\|w\| = \|v\|$ and $\|w - v\| \le \epsilon \|v\|$.
 \end{defn}

The following well-known theorem describes exactly which
connected, semisimple, real Lie groups have Kazhdan's
property~$(T)$. We note, in particular, that $\SL(2,\real)$
does not have Kazhdan's property~$(T)$.

\begin{thm}[{Kazhdan, Kostant, Serre, Wang}]
\label{SSgrps-propertyT}
 Let $G$ be a connected semisimple real Lie group.

 \begin{enumerate}
 \item \label{SSgrps-propertyT-simple}
 Assume $G$ is simple. Then $G$ has Kazhdan's
property~$(T)$ if and only if either 
 \begin{itemize}
 \item $\Rrank(G) \ge 2$ or
 \item  $G$ is compact, or 
 \item $G$ is locally isomorphic to either
$\operatorname{Sp}(1,n)$ or the real-rank one form of~$F_4$.
 \end{itemize}
 \item \label{SSgrps-propertyT-factors}
 $G$ has Kazhdan's property~$(T)$ if and only if each
simple factor of~$G$ has Kazhdan's property~$(T)$.
 \item \label{SSgrps-propertyT-center}
 $G$ has Kazhdan's property~$(T)$ if and only if
$G/Z(G)$ has Kazhdan's property~$(T)$.
 \end{enumerate}
 \end{thm}

\begin{pf}
 For~\pref{SSgrps-propertyT-simple}, see
\cite[Thm.~III.5.6(c), p.~131]{MargulisBook}. 
 For~(\ref{SSgrps-propertyT-factors}$\Leftarrow$), see
\cite[Cor.~III.2.10, p.~117]{MargulisBook}. 
 For~(\ref{SSgrps-propertyT-factors}$\Rightarrow$ and
\ref{SSgrps-propertyT-center}$\Rightarrow$), see
\cite[Lem.~III.2.4, p.~115]{MargulisBook}. 
 For (\ref{SSgrps-propertyT-center}$\Leftarrow$), see \cite[Thm.~2.12,
p.~28]{delaHarpeValette}.
 \end{pf}

In the proof of our generalization of the Thurston Stability
Theorem~\ref{C1-fp-measure}, the following lemma is used
to construct vectors~$v$ as in
Definition~\ref{KazhdanDefn}.

\begin{lem}[{cf.~\cite[2nd par.\ of pf.~of Thm.~9.1.1,
p.~163]{ZimmerBook}}] \label{localcontrol}
 Let $\alpha \colon G \times M \to \Diff^1(I)$ be a
Borel cocycle. For each $g \in G$, assume that for almost
every $m \in M$, we have $\alpha(g,m)(0) = 0$ and
$\alpha(g,m)'(0) = 1$.
 Then, for every compact subset~$C$ of~$G$, and
every $\epsilon > 0$, there is a nontrivial interval~$I'$
containing~$0$, such that, for every $g \in C$, we have
 $$ \mu \bigset{ m \in M}
 { \forall s \in I', ~ \bigl|\alpha(g,m)'(s) - 1 \bigr| <
\epsilon } > 1-\epsilon .$$
 \end{lem}

\section{A crucial lemma} \label{crucial-lemma}

Ghys' proof of Theorem~\ref{GhysC1LatticeThm} is based on
the existence certain subgroups $P$ and~$L$ of~$G$, such
that $P \subset L$, and the action of~$L$ on~$L/P$ is
triply transitive. (Then this is contrasted with the fact
that the group of orientation-preserving homeomorphisms
of~$\torus$ is not triply transitive on~$\torus$.) Ghys
describes $P$ and~$L$ quite explicitly, in geometric terms,
but this depends on a case-by-case study that uses the
classification of semisimple Lie groups. By giving a
uniform construction, the following lemma allows us to
avoid case-by-case analysis (or, at least, to condense it
into this one lemma). 

\begin{lem} \label{goodLs}
 Let 
 \begin{itemize}
 \item $H$ be a connected, noncompact, almost simple, real Lie
group;
 \item $P$ be a minimal parabolic subgroup of~$H$; and 
 \item $A$ be a maximal split torus of~$H$ contained
in~$P$. 
 \end{itemize}
 If $H$ is not locally isomorphic to $\SL(2,\real)$, then
there is a connected Lie subgroup~$L$ of~$H$, such that:
 \begin{enumerate}
 \item \label{goodLs-notinP}
 $L \not\subset P$;
 \item \label{goodLs-A}
 $\Lie A \cap [\Lie L,\Lie L]$ is nontrivial;
 \item \label{goodLs-C}
 $C_{\Lie A}(\Lie L)$ has codimension one in~$\Lie A$; and
 \item \label{goodLs-triply}
 $L N_P(L)$ is triply transitive on $L N_P(L)/N_P(L)$.
 \end{enumerate}
 \end{lem}

\begin{pf}
 Let us begin by making our goal more specific.

\begin{Claim}
 It suffices to find  a connected, closed subgroup~$L$
of~$H$, a real root~$\alpha$ of~$H$, and an element~$g$
of~$H$, such that:
 \begin{enumerate} \renewcommand{\theenumi}{\alph{enumi}}
 \item \label{goodLclaim-SL2}
 $L$ is locally isomorphic to $\SL(2,\real)$;
 \item \label{goodLclaim-alpha}
 $\Lie L = \langle \Lie H_\alpha \cap \Lie L , \Lie
H_{-\alpha} \cap \Lie L \rangle$; 
 \item \label{goodLclaim-centralize}
 $g \in C_H(A)$; and
 \item \label{goodLclaim-normalize}
 $g$ normalizes~$L$, and acts on~$L$ by an outer
automorphism.
 \end{enumerate}
 \end{Claim}

\begin{pf*}{Proof of Claim}
 \pref{goodLs-notinP} Because $P$ is minimal parabolic, we
know that $P/\Rad P$ is compact, so $P$ does not
contain~$L$ (or any other noncompact, semisimple subgroup).

\pref{goodLs-A} By definition, $[\Lie L,\Lie L]$ contains
the nontrivial subalgebra $[\Lie H_\alpha \cap \Lie L, \Lie
H_{-\alpha} \cap \Lie L]$ of~$\Lie A$.

\pref{goodLs-C} Because $\ker(\alpha)$ centralizes~$\Lie L$,
we know that $C_{\Lie A}(\Lie L)$ has codimension one
in~$\Lie A$.

 \pref{goodLs-triply} Because $P$ is parabolic and
contains~$A$, we know that $\Lie P$ contains either $\Lie
H_\alpha$ or~$\Lie H_{-\alpha}$ \see{P-structure}. Thus,
$\Lie P \cap \Lie L$ is a parabolic subalgebra of~$\Lie L$
\see{parab-sl2}. Then Lemma~\ref{SL2-almparab} implies that
$P \cap L$ is parabolic in~$L$. Thus, we may identify $L/(P
\cap L)$ with $\real P^1 \approx \torus$, so there is an
$L$-invariant circular order on $L/(P \cap L)$, and $L$ has
only two orbits on the ordered triples of distinct points
in $L/(P \cap L)$: the positively oriented triples and the
negatively oriented triples. Modulo inner automorphisms,
there is only one outer automorphism of~$L$, so it is easy
to verify that any outer automorphism of~$L$ that fixes $P
\cap L$ must take each positively oriented triple to a
negatively oriented triple. Thus,  because $g \in C_G(A)
\subset P$ \see{P-structure}, we see
from~\pref{goodLclaim-normalize} that all ordered triples
of distinct points in $L/(P \cap L)$ are in the same
$\bigl( L \, N_P(L) \bigr)$-orbit. 

This completes the proof of the claim.
 \renewcommand{\qed}{} 
 \end{pf*}

 We now consider two cases, based on the real rank of~$H$.

\begin{case}
 Assume $\Rrank H = 1$.
 \end{case}
 From the classification of simple Lie groups of real rank
one (cf.~\cite[Table~X.V, p.~518]{HelgasonBook}) (and the
fact that $H$ is not locally isomorphic to $\PSL(2,\real)
\iso \SO(1,2)$), we know that $H$ must contain a subgroup
locally isomorphic to either $\SO(1,3)$ (if $H$ is locally
isomorphic to $\SO(1,n)$) or $\SU(1,2)$ (if $H$ is locally
isomorphic to $\SU(1,n)$, $\Sp(1,n)$, or the rank one form
of~$F_4$). Then the proof is completed by explicitly
constructing $L$~and~$g$ for $\SO(1,3)$ and $\SU(1,2)$.

\setcounter{subcase}{0}

\begin{subcase}
 Assume $H$ is locally isomorphic to $\SO(1,3)$.
 \end{subcase}
 We may assume that $H = \SL(2,\complex)$, that $A$ consists
of diagonal matrices, and that $P$ is the group of upper
triangular matrices.
 The matrix
 $\left(\begin{smallmatrix}
 i&0\\
 0&-i\\
 \end{smallmatrix} \right)$
 acts by an outer automorphism of $\SL(2,\real)$.

\begin{subcase}
 Assume $H = \SU(1,2)$.
 \end{subcase}
 We use the Hermitian form
 $\langle x|y \rangle = x_1 \overline{y_3} + x_2
\overline{y_2} + x_3 \overline{y_1}$.
  We may assume that $A$ consists of diagonal matrices, and
that $P$ is the group of upper triangular matrices in~$H$.
 Let 
 $$\Lie L = 
 \bigset{ \begin{pmatrix}
 a&t&0\\
 s&1&-t\\
 0&-s&-a\\
 \end{pmatrix}
 }
 {a,s,t \in \real}
 \text{ and }
 g =\left( \begin{matrix}
 -1&0&0\\
 0&1&0\\
 0&0&-1\\
 \end{matrix} \right) .$$

\begin{case}
 Assume $\Rrank H > 1$.
 \end{case}
 It is well known (see, for example, \cite[Prop.~I.1.6.2,
p.~46]{MargulisBook}) that $H$ contains a closed,
connected subgroup that is locally isomorphic to either
$\SL(3,\real)$ or $\Sp(4,\real)$. Therefore, by passing to
a subgroup, and then  passing to a locally isomorphic
group, we may assume that $H$ is either $\SL(3,\real)$ or
$\Sp(4,\real)$.

\begin{subcase}
 Assume $H = \SL(3,\real)$.
 \end{subcase}
  We may assume that $A$ consists of diagonal matrices, and
that $P$ is the group of upper triangular matrices. Let
 $$ 
 L = 
 \begin{pmatrix}
 \SL(2,\real)&0 \\
 0&1 \\
  \end{pmatrix}
 .$$
 The matrix
 $g = \left(\begin{smallmatrix}
 1 &  0 & 0 \\
 0 & -1 & 0 \\
 0 &  0 & -1 \\
 \end{smallmatrix}\right)
 $
 acts by an outer automorphism of~$L$, and
centralizes~$A$.

\begin{subcase}
 Assume $H = \Sp(4,\real)$.
 \end{subcase}
 We use the symplectic form defined by
 $$ \langle (x_1,x_2) | (y_1,y_2) \rangle = x_1 \cdot y_2 -
x_2 \cdot y_1 ,$$
 for $x_i,y_i \in \real^2$, and we may assume that $A$
consists of diagonal matrices.
 Let
 $$L = 
 \bigset
 {\begin{pmatrix}
 R&0\\
 0& \theta(R)\\
 \end{pmatrix}
 }
 {R \in \SL(2,\real)} ,$$
 where $\theta$ is the Cartan involution (transpose-inverse).
  The matrix
 $g = \left(\begin{smallmatrix}
 1 &  0 & 0 & 0 \\
 0 & -1 & 0 & 0 \\
 0 &  0 & 1 & 0 \\
 0 &  0 & 0 &-1 \\
 \end{smallmatrix}\right)
 $
 acts by an outer automorphism of~$L$, and
centralizes~$A$.
 \end{pf}

\begin{rem} \label{goodLs-conceptual}
 By using more theory, one can give a more conceptual proof
of Lemma~\ref{goodLs}, without using the classification of
real simple Lie algebras.

\setcounter{case}{0} 

\begin{case} \label{L(rank1)}
 Assume $\Rrank H = 1$.
 \end{case}
 Write $P^\circ = CAU$, where $C$ is a compact, connected
subgroup of~$C_H(A)$ and $U$ is the unipotent radical
of~$P$. Let $\alpha$ be the simple real root of~$H$, and
assume without loss of generality that $\Lie H_\alpha
\subset \Lie U$. Because the compact, connected group~$C$
acts nontrivially on~$\Lie H_\alpha$, there is some $g \in
C$ and $u \in \Lie H_\alpha$, such that $\Ad g(u) = -u$.
From the Jacobson-Morosov Theorem, we know that $u$ is
contained in a subalgebra~$\Lie L$ that is isomorphic to
$\Sl(2,\real)$.

Since $\Rrank H = 1$, we know that $N_{\Lie H}(\langle u
\rangle) \subset \Lie P$, so $\Lie P \cap \Lie L$ contains
a maximal split torus of~$\Lie L$. Thus, because all
maximal split toral subalgebras of~$\Lie P$ are conjugate,
there is some $v \in U$, such that $(\Ad v)(\Lie A)$ is a
maximal split toral subalgebra of~$\Lie L$ that
normalizes~$\langle u\rangle$. Then $\Lie A$ normalizes
$(\Ad v^{-1})\langle u \rangle$, so $(\Ad v^{-1})u \in \Lie
H_\alpha$. Because $u$ is also in~$\Lie H_\alpha$, and
$[\Lie U, \Lie U] \cap \Lie H_\alpha = 0$, we conclude that
$(\Ad v^{-1})u = u$. Thus, replacing $\Lie L$ by $(\Ad
v^{-1})\Lie L$, we may assume that $\Lie A \subset \Lie L$. 

Then $g$ normalizes the parabolic subalgebra $\Lie A +
\langle u \rangle$ of~$\Lie L$, so it must normalize~$\Lie
L$. Also, we know that $g$ acts on~$\Lie L$ by an outer
automorphism, because $g$ conjugates~$u$ to~$-u$, whereas no
nontrivial unipotent element is conjugate to its inverse in
$\SL(2,\real)$.

\begin{case} \label{L(rank2)}
 Assume $\Rrank H > 1$.
 \end{case}
 For simplicity, let us assume that $H$ is $\real$-split.
Choose two roots $\alpha$~and~$\beta$, such that the
$\beta$-string through~$\alpha$ has odd length, let $\Lie L
= \langle \Lie H_\alpha, \Lie H_{-\alpha} \rangle$, let
$L_\beta$ be the connected Lie subgroup of~$H$ corresponding
to the subalgebra $\langle \Lie H_\beta, \Lie H_{-\beta}
\rangle$,  and let $V$ be the $L_\beta$-submodule of~$\Lie
H$ generated by~$\Lie H_\alpha$. Then, identifying
$L_\beta$ with $\SL(2,\real)$, the highest weight of~$V$ is
odd, so $g = 
 \begin{pmatrix} -1 & 0 \\ 0 & -1 \end{pmatrix}$
 acts as~$-1$ on the highest weight space~$\Lie H_\alpha$. \qed
 \end{rem}

\section{Proof of Theorem~\ref{GhysThm}} \label{ProveGhys}

The reader is encouraged to read Ghys' beautiful proof
\cite[\S4]{Ghys} for the case of lattices in $\SL(3,\real)$
before looking at the general case considered here.
 Many of the ideas of this section can be found
in~\cite{Ghys}, but we have reorganized them, and changed
some of the emphasis. Ghys' proof is presented in geometric
terms, but we have reformulated the argument in
group-theoretic terms.

\begin{notation} \mbox{ }
 \begin{itemize}
 \item $G$, $M$, and~$\alpha$ are always assumed to be as
described in the statement of Theorem~\ref{GhysThm}. (In
particular, $G$ has no factors locally isomorphic to $\SL(2,
\real)$.)
 \item $P$ is a minimal parabolic subgroup of~$G$.
 \item For any natural number~$k$, $\torus_k$ denotes the
collection of all $k$-element subsets of~$\torus$.
 \end{itemize}
 \end{notation}

\begin{lem} \label{GhysPf-orpres}
 We may assume $\alpha \colon G \times M \to
\Homeo_+(\torus)$.
 \end{lem}

\begin{pf}
 Let $\sgn \colon \Homeo(\torus) \to \{\pm1\}$ be the
homomorphism with kernel $\Homeo_+(\torus)$, and let
$\varepsilon = \sgn \circ \alpha$, so $\varepsilon \colon G
\times M \to \{\pm1\}$ is a Borel cocycle.

Let $M^+ = M \mathbin{\times_\varepsilon} \{\pm1\}$.
Because $M^+$ is a two-point extension of~$M$ and $M$ is
irreducible, it is clear that each closed,
connected, noncompact, normal subgroup of~$G$ has no more
than two ergodic components on~$M^+$. We may assume that
$G$ is ergodic on~$M^+$, for, otherwise, $\varepsilon$ is
equivalent to the trivial cocycle, so $\alpha$ is
equivalent to a cocycle into $\Homeo_+(\torus)$, as
desired. Then $G$ must act ergodically on the space of
ergodic components of any normal subgroup. Because $G$,
being connected, has no nontrivial action on any finite
set, we conclude that $M^+$ is irreducible. 

Define $\alpha^+ \colon G \times M^+ \to \Homeo(\torus)$
by $\alpha^+ \bigl( g, (m, \pm 1) \bigr) = \alpha(g,m)$.
If there is a $G$-invariant probability measure~$\nu^+$ on
$M^+ \mathbin{\times_{\alpha^+}} \torus$, such that $\nu^+$
projects to~$\mu^+$ on~$M^+$, then simply let $\nu$ be the
projection of~$\nu^+$ to $M \atimes \torus$.

Now let $f$ be any
orientation-reversing homeomorphism of~$\torus$, and 
 define $\sigma \colon M^+ \to \Homeo(\torus)$ by 
 $$ \sigma ( m, \varepsilon ) = 
 \begin{cases}
 \Id & \text{ if $\varepsilon = 1$} \\
 f & \text{ if $\varepsilon = -1$} 
 \end{cases} .$$
 For any $m^+ \in M^+$, we have
 $\sigma(g m^+) \alpha^+(g, m^+) \sigma(m^+)^{-1} \in
\Homeo_+(\torus)$, so we see that $\alpha^+$
is cohomologous (via~$\sigma$) to a cocycle with values in
$\Homeo_+(\torus)$.
 \end{pf}

Henceforth, we assume $\alpha(G \times M) \subset
\Homeo_+(\torus)$.

It suffices to show that there is an $\alpha$-equivariant
Borel map $\psi \colon M \to \Prob(\torus)$, for then we
may set
 $\nu = \int_M \bigl( m \times \psi(m) \bigr) \, d\mu(m)$.
 The action of~$G$ on~$(G/P) \times M$ is amenable (because
$P$ is amenable) \cite[4.1.7bis, 4.3.2, 4.3.4]{ZimmerBook},
and the space of measurable functions from $(G/P) \times M$
to $\Prob(\torus)$ is an affine $G$-space over $(G/P)
\times M$ \cite[Defn~4.3.1]{ZimmerBook}. Thus, from the
definition of an amenable action
\cite[Defn~4.3.1]{ZimmerBook}, we know that there is an
$\alpha$-equivariant Borel map $\Psi \colon (G/P) \times M
\to \Prob(\torus)$ (cf.~\cite[pf.\  of Step~1 of
Thm.~5.2.5, bot.\ of p.~103]{ZimmerBook}). The following
theorem completes the proof of Theorem~\ref{GhysThm}.

\begin{thm} \label{Ghys-essconst}
 Suppose $\Psi \colon (G/P) \times M \to \Prob(\torus)$ is
an $\alpha$-equivariant Borel map. For each $m \in M$,
define $\Psi_m \colon G/P \to \Prob(\torus)$ by $\Psi_m(x)
= \Psi(x,m)$.

 Then $\Psi_m$ is essentially constant, for a.e.~$m \in M$.
 \end{thm}

\begin{pf}
 If almost every $\Psi(x,m)$ is atomless, the desired
conclusion is given by Theorem~\ref{Ghys-atomless} below. If
there is some~$k$, such that almost every $\Psi(x,m)$
consists of $k$~atoms of equal weight, the desired conclusion
is given by Corollary~\ref{Ghys-atoms} below. Because $G$ is
ergodic on $(G/P) \times M$ \see{MooreErgodicityThm}, it is
not difficult to reduce the problem to these two cases.

Namely, any $\nu \in \Prob(\torus)$ has a unique
decomposition of the form $\nu = \nu_0 + \nu_1$, where
$\nu_0$ has no atoms, and $\nu_1$ consists entirely of
atoms. (Either of the terms in the decomposition may
be~$0$.) Thus, we may write $\Psi = \Psi_0 + \Psi_1$, where
$\Psi_i(x,m) = [\Psi(x,m)]_i$. Because the decomposition
$\nu = \nu_0 + \nu_1$ is $\Homeo(\torus)$-equivariant (and
unique), we see that $\Psi_0$ and~$\Psi_1$ are
$\alpha$-equivariant. Then, because $G$ is ergodic on
$(G/P) \times M$, we see, for $i = 0,1$, that either
$\Psi_i = 0$ for a.e.~$(x,m)$ or $\Psi_i \neq 0$ for
a.e.~$(x,m)$. Thus, either $\Psi_i = 0$~a.e.\ (in which
case $\Psi = \Psi_{1-i}$), or, after renormalizing, $\Psi_i$
defines an $\alpha$-equivariant Borel map into
$\Prob(\torus)$. Then, because the sum of
$\alpha$-equivariant functions is $\alpha$-equivariant,
there is no harm in assuming that either $\Psi = \Psi_0$ or
$\Psi = \Psi_1$.

If $\Psi = \Psi_0$, then Theorem~\ref{Ghys-atomless}
shows that $\Psi_m$ is essentially constant.

Thus, we henceforth assume that $\Psi = \Psi_1$.
For any $\nu \in \Prob(\torus)$ that consists entirely of
atoms, and any rational number $q \in (0,1)$, let $\nu^{>q}
\subset \torus$ be the set of atoms of weight~$>q$. Because
this definition is $\Homeo(\torus)$-equivariant, and $G$ is
ergodic on $(G/P) \times M$, we see that the cardinality
of~$\Psi^{>q}$ is constant~a.e., so $\Psi^{>q}$ is an
$\alpha$-equivariant Borel map into~$\torus_k$, for
some~$k$. Then Corollary~\ref{Ghys-atoms} asserts that
$\Psi_m^{>q}$ is essentially constant. Because this is true
for all rational~$q$, we conclude that $\Psi_m$ itself is
essentially constant, as desired.
 \end{pf}

\begin{thm} \label{Ghys-atomless}
 Suppose $\Psi \colon (G/P) \times M \to \Prob(\torus)$ is
an $\alpha$-equivariant Borel map. For each $m \in M$,
define $\Psi_m \colon G/P \to \Prob(\torus)$ by $\Psi_m(x)
= \Psi(x,m)$.

If $\Psi(x,m)$ is atomless, for almost every $(x,m) \in
(G/P) \times M$, then $\Psi_m$ is essentially constant, for
a.e.~$m \in M$.
 \end{thm}

\begin{pf}
 Let $\Prob_0(\torus)$ be the set of atomless probability
measures on~$\torus$.
 Define
 $$ \Psi^2 \colon (G/P)^2 \times M \to \Prob_0(\torus)^2
 \text{ \ by \ }
 \Psi^2(x_1,x_2,m) = \bigl( \Psi(x_1,m), \Psi(x_2,m) \bigr)
$$
 and
 $$ D \colon \Prob_0(\torus)^2 \to [0,1]
 \text{ \ by \ }
 D(\mu_1,\mu_2) = \sup_J \bigl| \mu_1(J)  - \mu_2(J)
 \bigr| ,$$
 where $J$ ranges over all subintervals of~$\torus$.
 It suffices to show that the composite function
 $D \circ \Psi^2 \colon (G/P)^2 \times M \to [0,1]$
 is~$0$~a.e.

\setcounter{step}{0}

\begin{step}
 $D$ is continuous.
 \end{step}
 Given $\mu_1,\mu_2 \in \Prob_0(\torus)$. Because $\mu_1$
and~$\mu_2$ are atomless, there is a mesh
$t_0,t_1,\ldots,t_n=t_0$ of points in~$\torus$, such that
$\mu_k([t_i,t_{i+1}])< \epsilon/40$, for each~$i$ and for
$k = 1,2$. Also, for each $i,j \in \{0,\ldots,n\}$, there
are continuous functions $f^+_{ij},f^-_{ij} \colon
\torus \to [0,1]$, such that 
 $\supp f^-_{ij} \subset (t_i,t_j)$,
 $f^+_{ij} \bigl( [t_i,t_j] \bigr) = 1$, and
 $\mu_k(f^+_{ij} - f^-_{ij}) < \epsilon/40$ for $k = 1,2$.
 If
$\nu_k$ is a measure so close to~$\mu_k$ that
 $|\nu_k(f^\varepsilon_{ij}) - \mu_k(f^\varepsilon_{ij})| <
\epsilon/40$ for all $i,j \in \{0,\ldots,n\}$ and
$\varepsilon \in \{+,-\}$, then 
 $$ \Bigl| \nu_k \bigl( [t_i,t_j] \bigr) - \mu_k
\bigl( [t_i,t_j] \bigr) \Bigr| < \frac{\epsilon}{20}$$
 for all $i$~and~$j$. Therefore 
 $\bigl| \nu_k(J)  - \mu_k(J)
 \bigr| < \epsilon/2$
 for every interval~$J$, so
 $|D(\nu_1,\nu_2) - D(\mu_1,\mu_2) | \le \epsilon$. This
proves the continuity of~$D$.

\begin{step}
 $D \circ \Psi^2$ is essentially constant.
 \end{step}
 Because $\Psi^2$ is $\alpha$-equivariant and $D$ is
$\Homeo(\torus)$-invariant, we know that $D \circ \Psi^2$ is
essentially $G$-invariant.  The Moore Ergodicity
Theorem~\ref{MooreCor(G/P)^2} implies that $G$ is
ergodic on $(G/P)^2 \times M$, so we conclude that $D \circ
\Psi^2$ is essentially constant. 

\begin{step}
 We have $D \circ \Psi^2 = 0$ a.e.
 \end{step}
From Lusin's Theorem, we know that
$\Psi$ is continuous on some compact subset~$C$ of positive
measure in $G/P$. Therefore, $D \circ \Psi^2$ is continuous
on $C \times C$. By replacing~$C$ with a smaller compact
set, we may assume that every conull subset of~$C$ is
dense. Then, because $D \circ \Psi^2$ is essentially
constant, we conclude that $D \circ \Psi^2$ is 
constant on $C \times C$. Obviously, $D \circ \Psi^2$ is~$0$
on the diagonal $\{(c,c)\}$, so we conclude that  $D \circ
\Psi^2$ is~$0$~a.e.
 \end{pf}

\begin{thm} \label{SL2inv}
 Suppose $\Psi \colon (G/P) \times M \to \torus_k$ is an
$\alpha$-equivariant Borel map. For each $m \in M$, define $\Psi'_m \colon G \to \torus_k$
by $\Psi'_m(g) = \Psi(gP,m)$.

 Suppose $L$ is a closed, connected subgroup of~$G$, such
that
 \begin{enumerate}
 \item \label{SL2inv-noncpct}
 $C_P(L)$ is not compact; and
 \item \label{SL2inv-triply}
 $L \, N_P(L)$ acts triply transitively on~$L \, N_P(L) /
N_P(L)$.
 \end{enumerate}
 Then $\Psi'_m$ is essentially right $L$-invariant, for
a.e.~$m \in M$. {\upshape(}That is,
for each $l \in L$, we have $\Psi'_m(g l) = \Psi'_m(g)$ for
a.e.~$g \in G$.{\upshape)}
 \end{thm}

\begin{pf}
 Because $\Psi'_m$ is right $P$-invariant, we may assume that
$L \not\subset P$.

 The inclusion $N_P(L) \hookrightarrow L \, N_P(L)$ induces
a $G$-equivariant smooth submersion
 $$ \pi \colon G/ N_P(L) \to G/ \bigl( L \, N_P(L) \bigr)
.$$
 Define
 $$ \overline{X} =
 \bigset{ (x_1,x_2,x_3) \in \bigl( G/ N_P(L) \bigr)^3 }
 { \pi(x_1) = \pi(x_2) = \pi(x_3) }$$
 and
 $$ X =
 \bigl\{\, (x_1,x_2,x_3) \in \overline{X} \mid
 \mbox{$x_1, x_2, x_3$ distinct} \,\bigr\}
 .$$
 Then $\overline{X}$ is a closed submanifold of $\bigl( G/ N_P(L)
\bigr)^3$, and $X$ is conull open subset of~$\overline{X}$ (with
respect to any smooth measure on~$\overline{X}$).
 For $i=1,2,3$, let $\pi_i \colon \overline{X} \to G/P$ be the
$G$-equivariant map defined by $\pi_i(x_1,x_2,x_3) = x_i P$.
 
Because $G$ is transitive on $G/N_P(L)$, any $G$-orbit
on~$X$ contains a point $(x_1,x_2,x_3)$, such that $x_1 =
N_P(L)$. Then $x_1,x_2,x_3$ are three points in
 $L \, N_P(L) / N_P(L)$.
 Thus, assumption~\pref{SL2inv-triply} implies that $G$ is
transitive on~$X$. In particular, this implies that the
class~$\chi$ of Lebesgue measure is the unique
$\sigma$-finite $G$-invariant measure class on~$X$. For $i
= 1,2,3$, the projection $(\pi_i)_* \chi$ must be the
$G$-invariant measure class on $G/P$. Thus, we have an
essentially well-defined Borel map
 $\Psi^3 \colon \overline{X} \times M \to (\torus_k)^3$
given by
 $$\Psi^3(x,m) = \Bigl( \Psi\bigl( \pi_1(x),m \bigr),
\Psi\bigl( \pi_2(x),m \bigr), \Psi\bigl( \pi_3(x),m \bigr)
\Bigr) .$$
 Note that $\Psi^3$ is $\alpha$-equivariant.

The stabilizer of a triple of points in $L \, N_P(L)/
N_P(L)$ obviously contains $C_P(L)$, which,
by~\pref{SL2inv-noncpct}, is not compact. Thus, we conclude
from the Moore Ergodicity Theorem~\ref{MooreErgodicityThm}
that $G$ is ergodic on $\overline{X} \times M$. This
implies that there is a single $\Homeo(\torus)$-orbit~$O$
on~$(\torus_k)^3$, such that $\Psi^3(x,m) \in O$ for
a.e.~$(x,m)$.
 For any permutation~$\sigma$ of~$\{1,2,3\}$, and any
$(x_1,x_2,x_3) \in \overline{X}$, we know that
$(x_{\sigma(1)},x_{\sigma(2)},x_{\sigma(3)})$ also belongs
to~$\overline{X}$. Therefore,  Lemma~\ref{circleorder}
implies that 
 $O = \{ \, (A,A,A) \mid A \in \torus_k \,\}$.

The map $G \times L^2 \to \overline{X}$ given by
 $(g,l,l') \mapsto \bigl( g N_P(L),  gl N_P(L),  gl' N_P(L)
\bigr)$ is a submersion, so it preserves the class of
Lebesgue measure. Thus, from the conclusion of the preceding
paragraph, we see that, for almost every $m \in M$, $g \in
G$, and $l,l' \in L$, we have $\Psi'_m(g) = \Psi'_m(gl) =
\Psi'_m(gl')$. From Fubini's Theorem (and ignoring~$l'$), we
conclude, for a.e.~$m \in M$, that $\Psi'_m$ is essentially
right $L$-invariant.
 \end{pf}

\begin{cor} \label{Ghys-atoms}
 Suppose $\Psi \colon (G/P) \times M \to \torus_k$ is an
$\alpha$-equivariant Borel map.
 For each $m \in M$, define $\Psi_m \colon G/P \to \torus_k$
by $\Psi_m(x) = \Psi(x,m)$.

 Then $\Psi_m$ is essentially constant, for a.e.~$m \in M$.
 \end{cor}

\begin{pf}
 Let $P = MAN$ be the Langlands decomposition of~$P$
\cite[p.~81]{WarnerBook}. (Because the parabolic
subgroup~$P$ is minimal, we know that $A$ is a maximal
split torus of~$G$.)

 It suffices to show, for each simple factor~$H$ of~$G$,
that there are subgroups $L_1,L_2,\ldots,L_n$ of~$H$, such
that 
 \begin{enumerate}
 \renewcommand{\theenumi}{\alph{enumi}} 
 \item \label{goodLi}
 each subgroup~$L_i$ satisfies the hypotheses of
Theorem~\ref{SL2inv}, and
 \item \label{LcapAgenerates}
 $\{[L_i,L_i] \cap A\}$ generates~$A \cap H$.
 \end{enumerate}
 To see that this suffices, let $J$ be the subgroup generated
by $\{P \cap H\} \cup \{L_1, L_2,\ldots,L_n\}$. Then
Theorem~\ref{SL2inv} implies that $\Psi'_m$ is essentially
right $J$-invariant. Because $J \supset P \cap H$, we know
that $J$ is parabolic in~$H$; let $J
= M_J A_J N_J$ be the Langlands decomposition of~$J$, with
$A_J \subset A \cap H$. Then $[J^\circ,J^\circ] \subset M_J
N_J$, so $[J^\circ, J^\circ] \cap A_J = e$. 
On the other hand, we have $(A \cap H)^\circ \subset
[J^\circ,J^\circ]$ (see~\pref{LcapAgenerates}).
We conclude that $A_J$ is trivial, so 
 $$ J \supset M_J A_J = C_H(A_J) =
C_H(e) = H .$$
 Therefore $\Psi'_m$ is essentially right
$H$-invariant. Because this is true for every simple
factor~$H$, we conclude that $\Psi'_m$ is essentially right
$G$-invariant, so $\Psi'_m$ is essentially constant, for
a.e.~$m \in M$.

If $\Rrank H = 1$, then Lemma~\ref{goodLs} provides an
appropriate subgroup~$L$ satisfying \pref{goodLi}
and~\pref{LcapAgenerates}. (Because $L$ is centralized by
all the simple factors other than~$H$, the requirement that
$C_P(L)$ be noncompact is automatically satisfied.)

We may now assume that $\Rrank H > 1$. Let $L$ be as in
Lemma~\ref{goodLs}, and let $W$ be the relative Weyl group
of~$\Lie H$ (with respect to~$\Lie A \cap \Lie H$). Because
$[\Lie L,\Lie L] \cap \Lie A$ is nontrivial and $W$ acts
irreducibly on~$\Lie A$, we know that
 $\bigset{w \bigl( [\Lie L,\Lie L] \cap \Lie A \bigr)}{w
\in W}$
 spans~$\Lie A$, so $\{\, w(L) \mid w \in W \,\}$
satisfies~\pref{LcapAgenerates}. Because $C_{\Lie A} \bigl(
w(\Lie L) \bigr)$ has codimension one in~$\Lie A$ (and,
being a subspace of~$\Lie A$, is contained in~$\Lie P$), we
know that $C_P(L)$ is noncompact. Thus, we see that each
$w(L)$ satisfies the hypotheses of Theorem~\ref{SL2inv}.
 \end{pf}

 The following result was used in the proof of
Theorem~\ref{SL2inv}. For completeness, we include
the proof. We also remark that, as explained by
Ghys~\cite[Step~3 of~\S4, bot.\ of p.~210]{Ghys}, the group
$\Homeo_+(\torus)$ has only finitely many orbits
on~$(\torus_k)^3$.

\begin{lem}[{Ghys}] \label{circleorder}
 Let $O$ be an orbit of $\Homeo_+(\torus)$
on~$(\torus_k)^3$, and assume there is an element
$(A_1,A_2,A_3)$ of~$O$, such that $\left(
A_{\sigma(1)},A_{\sigma(2)},A_{\sigma(3)} \right) \in O$,
for every permutation~$\sigma$ of $\{1,2,3\}$. Then 
 $O = \{\, (A,A,A) \mid A \in \torus_k \,\}$.
 \end{lem}

\begin{pf*}{Proof {\rm \cite[bot.\ of p.~211]{Ghys}}}
 Let $B = A_1 \cup A_2 \cup A_3$, and let 
 $H = \{\, h \in \Homeo_+(\torus) \mid h(B) = B \,\}$. For
each permutation~$\sigma$ of $\{1,2,3\}$, there is an
orientation-preserving homeomorphism~$h_\sigma$ (not
unique) of~$\torus$, such that $h_\sigma(A_i) =
A_{\sigma(i)}$. Then $h_\sigma \in H$, and the restriction
of~$H$ to~$B$ is a cyclic group, so the commutator of any
two of these homeomorphisms acts trivially on~$B$. Because
the permutation $\sigma = (1,2,3)$ is a commutator in the
symmetric group~$S_3$, we conclude that $h_{(1,2,3)}$ acts
trivially on~$B$. Because $h_{(1,2,3)}(A_1) = A_2$ and
$h_{(1,2,3)}(A_2) = A_3$, this implies $A_1 = A_2 = A_3$.
 \end{pf*}

\section{The Reeb-Thurston Stability Theorem}
\label{Thurston-section}

Ghys' proof of Theorem~\ref{GhysC1LatticeThm} relies on the
following one-dimensional case of the Reeb-Thurston
Stability Theorem~\cite{Thurston}. (See
\cite{ReebSchweitzer} and \cite{Schachermayer} for elegant
proofs.) 

\begin{thm}[{Thurston \cite{Thurston}}] \label{ReebThurston}
 If $\Gamma$ is a finitely generated group, such that
$\Gamma/[\Gamma,\Gamma]$ is finite, then there is no
nontrivial homomorphism $\Gamma \to \Diff^1_+(I)$.
 \end{thm}

For the proof of Theorem~\ref{C1-isometry}, we provide the
following generalization in the setting of Borel
cocycles. Applying this result to $G/\Gamma$ recovers
Thurston's theorem in the special case where $\Gamma$ is a
lattice in~$G$, and $G$ has Kazhdan's property~$(T)$
\see{KazhdanDefn}.

\begin{thm} \label{Kazhdan-homeo}
 Let 
 \begin{itemize}
 \item $G$ be a locally compact group with Kazhdan's
property~$(T)$;
 \item $M$ be a Borel $G$-space with invariant probability
measure~$\mu$; and
 \item $\alpha \colon G \times M \to \Diff^1_+(I)$ be a
Borel cocycle.
 \end{itemize}
 Then there is a $G$-invariant probability measure~$\nu$ on
$M \atimes I$, such that $\nu$ is equivalent to $\mu \times
\lambda$ {\upshape(}and $\nu$ projects to~$\mu$
on~$M${\upshape)}.

 Therefore, as a cocycle into $\HomeoLeb_+(I)$, $\alpha$ is
cohomologous to the trivial cocycle.
 \end{thm}

Before proving Theorem~\ref{Kazhdan-homeo}, let us explain
how it implies Theorem~\ref{C1-isometry}.

\begin{cor} \label{C1-fp-measure}
 Let 
 \begin{itemize}
 \item $G$ be a locally compact group with Kazhdan's
property~$(T)$;
 \item $M$ be a Borel $G$-space with finite invariant
measure~$\mu$; 
 \item $\alpha \colon G \times M \to \Diff^1_+(\torus)$ be a
Borel cocycle; and
 \item $f \colon M \to \torus$ be an
$\alpha$-equivariant measurable map.
 \end{itemize}
 Then there is a $G$-invariant probability measure~$\nu$ on
$M \atimes \torus$, such that $\nu$ is equivalent to $\mu
\times \lambda$.

 Therefore, as a cocycle into $\HomeoLeb_+(\torus)$,
$\alpha$ is cohomologous to the trivial cocycle.
 \end{cor}

\begin{pf}
 Cutting $\torus$ open at the point $f(m)$ yields an
interval~$I_m$, so we may define a cocycle $\hat\alpha
\colon G \times M \to \Diff^1_+(I)$. Then
Theorem~\ref{Kazhdan-homeo} applies.
 \end{pf}

\begin{thm} \label{C1-isometry}
  Let 
 \begin{itemize}
 \item $G$ be a connected, semisimple, real Lie group, such that
 \begin{itemize}
 \item $G$ has Kazhdan's property~$(T)$, and 
 \item $\Rrank G \ge 2$;
 \end{itemize}
 \item $M$ be an irreducible ergodic $G$-space with finite
invariant measure~$\mu$; and
 \item $\alpha \colon G \times M \to \Diff^1_+(\torus)$ be a
Borel cocycle.
 \end{itemize}
 Then there is a $G$-invariant probability measure~$\nu$ on
$M \atimes \torus$, such that $\nu$ is equivalent to $\mu
\times \lambda$.

 Therefore, as a cocycle into $\HomeoLeb_+(\torus)$,
$\alpha$ is cohomologous to a cocycle with values in the
rotation group $\Rot(\torus)$.
 \end{thm}

\begin{pf}
 Because $G$ has Kazhdan's property~$(T)$, we know that $G$
has no factors locally isomorphic to $\SL(2,\real)$
\see{SSgrps-propertyT}. Therefore, Theorem~\ref{GhysThm}
implies that there is a $G$-invariant probability
measure~$\sigma$ on $M \atimes \torus$, such that $\sigma$
projects to~$\mu$ on~$M$.

Define a cocycle $\beta \colon G \times (M \atimes \torus)
\to \Diff^1_+(\torus)$ by $\beta(g,m,s) = \alpha(g,m)$. The
map $f \colon M \atimes \torus \to \torus$ defined by
$f(m,s) = s$ is $\beta$-equivariant, so we know, from
Corollary~\ref{C1-fp-measure}, that there is a
$G$-invariant probability measure~$\hat\nu$ on $(M \atimes
\torus) \btimes \torus$, such that $\hat\nu$ is equivalent
to $\sigma \times \lambda$. 

Let $\nu$ be the image of $\hat\nu$ under the projection
$(m,s,t) \mapsto (m,t)$. Because $\hat\nu$ is equivalent to
$\sigma \times \lambda$, and $\sigma$ projects to~$\mu$
on~$M$, we see that $\nu$ is equivalent to $\mu \times
\lambda$, as desired.
 \end{pf}

To motivate the proof of Theorem~\ref{Kazhdan-homeo}, let
us sketch the analogous proof of
Theorem~\ref{ReebThurston}, under the assumption that
$\Gamma$ has Kazhdan's property~$(T)$. (It is well known
that, because $\Gamma$ is discrete, Kazhdan's
property~$(T)$ implies both that $\Gamma$ is finitely
generated \cite[Thm.~7.1.5, p.~131]{ZimmerBook} and that
$\Gamma/[\Gamma,\Gamma]$ is finite \cite[Cor.~7.1.7,
p.~131]{ZimmerBook}.)

\begin{pf*}{Proof of Theorem~\ref{ReebThurston} when
$\Gamma$ has Kazhdan's property~$(T)$}
 It suffices to show that the set of fixed points
of~$\Gamma$ is dense in~$I$. Suppose not. Then, replacing
$I$ by the closure of a component of the complement of the
fixed-point set, we may assume that there are no fixed
points in the interior of~$I$.

 We have a unitary representation~$\rho$ of~$\Gamma$ on
$L^2(I)$ given by
 $$ \bigl( \rho(\gamma)f \bigr) (t) = [\gamma'(t)]^{1/2}
f(\gamma^{-1} t) .$$

Let $\epsilon = 1/2$, and let $C \subset \Gamma$ and
$\delta > 0$ be as in Definition~\ref{KazhdanDefn}.
 Because $\Gamma/[\Gamma,\Gamma]$ is finite, the
homomorphism $\Gamma \to \real^+ \colon \gamma \mapsto
\gamma'(0)$ must be trivial. Thus, $\gamma'(0) = 1$, for
every $\gamma \in \Gamma$. Therefore, there is a nontrivial
interval~$I'$ containing~$0$, such that $|\gamma'(t) - 1| <
\delta^2/4$, for every $\gamma \in C$ and every $t \in I'$.
Let $\chi$ be the characteristic function of~$I'$. Then $\|
\rho(\gamma) f - f \| < \delta \|f\|$, for every $\gamma
\in C$, so we conclude from the choice of $C$ and~$\delta$
that there is some nonzero $\rho(\Gamma)$-invariant
function~$\phi$ in $L^2(I)$. Then $|\phi|^2 \, d\lambda$ is
a $\Gamma$-invariant measure on~$I$, so every point in the
support of this measure is fixed by~$\Gamma$. This
contradicts the assumption that $\Gamma$ has no fixed
points in the interior of~$I$.
 \end{pf*}

Our proof of Theorem~\ref{Kazhdan-homeo} is a fairly
straightforward translation of this argument to the setting
of cocycles for Borel actions, except that it is not
convenient to use a topological argument in this setting.
Therefore, instead of obtaining a contradiction
by finding a fixed point that does not belong to
the closure of the fixed-point set, we find a set
of fixed points whose measure is greater than the
measure of the set of fixed points.

\begin{pf*}{{\bf Proof of Theorem~\ref{Kazhdan-homeo}}}
 By passing to ergodic components, we may assume that $G$
is ergodic on~$M$.

 The map $G \times M \to \real^+$ defined by $(g,m) \mapsto
\alpha(g,m)'(0)$ is a cocycle. Because $G$ has Kazhdan's
property~$(T)$, and $\real^+$ is amenable and has no
compact subgroups, the cocycle must be cohomologous to the
trivial cocycle \cite[Thm.~9.1.1, p.~162]{ZimmerBook}, so,
by replacing $\alpha$ with an equivalent cocycle, we may
assume, for each $g \in G$, that $\alpha(g,m)'(0) = 1$ for
a.e.~$m \in M$.

 Because $\mu$ is $G$-invariant, we have
 $$ \int_{M \atimes I} \psi \, d(\mu \times \lambda)
 =  \int_{M \atimes I} \psi\bigl( g(m,s) \bigr)
  \, \alpha(g,m)'(s) \,\, d(\mu \times \lambda)(m,s) $$
 for any $g \in G$ and $\psi \in L^1(M \atimes I)$.
 Therefore, a unitary representation~$\rho$ of~$G$
on $L^2(M \atimes I)$ is given by
 $$ \bigl( \rho(g) \phi \bigr)(m,s) = 
 \phi \bigl( g^{-1}(m,s) \bigr)
 \, \bigl( \alpha(g^{-1},m)'(s) \bigr)^{1/2} $$
 for $g \in G$, $\phi \in L^2(M \atimes I)$, and
$(m,s) \in M \atimes I$.

 Fix a compact subset~$C$ of~$G$, as in
Definition~\ref{KazhdanDefn}.
 \begin{itemize}
 \item Fix some $\epsilon_1 \in (0,1)$, and let $\delta_1 =
\delta(\epsilon_1) > 0$ be the corresponding $\delta$-value
given by Definition~\ref{KazhdanDefn}. 
 \item Choose $\epsilon_2 > 0$ small enough that $9
\epsilon_2 < \delta_1^2$, and let $\delta_2 =
\delta(\epsilon_2) > 0$ be the corresponding $\delta$-value
given by Definition~\ref{KazhdanDefn}.
 \item Choose $\epsilon_3 > 0$ small enough that $13
\epsilon_3 < \delta_2^2$.
 \end{itemize}
 We may assume that $1 > \epsilon_1 > \delta_1 > \epsilon_2
> \delta_2 > \epsilon_3 > 0$.

Lemma~\ref{localcontrol} tells us that there is a
nontrivial interval~$I'$ containing~$0$, such that, for
every $g \in C$, we have
 \begin{equation} \equationlabel{mu(alpha-near-1)}
 \mu \bigl( \{ m \in M \mid \forall s \in
2I', ~ 
 |\alpha(g^{-1},m)'(s) - 1| < \epsilon_3 \} \bigr) > 1 -
\epsilon_3 .
 \end{equation}

Let $\invariant$ be the space of all $\rho(G)$-invariant
functions in $L^2(M \atimes I)$, and choose $\phi \in
\invariant$, such that
 $(\mu \times \lambda) \bigl( \phi^{-1}(0) \bigr)$ is
minimal.
 The minimum exists, because any convex combination of (the
absolute values of) countably many $\rho(G)$-invariant
functions is $\rho(G)$-invariant. Furthermore, 
 \begin{equation} \equationlabel{phi-1(0)}
 \text{for every $\psi \in \invariant$, we have $\psi = 0$
a.e.~on $\phi^{-1}(0)$.}
 \end{equation}

 Because $\phi$ is $\rho(G)$-invariant, we know that $\nu =
|\phi|^2 \cdot (\mu \times \lambda)$ is a $G$-invariant
measure on $M \atimes I$. (\emph{A priori,} $\phi$ could be
identically~$0$, so this measure could be trivial.) To
complete the proof, we will show that this measure is
equivalent to $\mu \times \lambda$; that is, we will show
that $(\mu \times \lambda) \bigl( \phi^{-1}(0) \bigr) = 0$. 
(Then $\nu$ projects to~$\mu$ on~$M$. Indeed, because $G$
is ergodic on $(M,\mu)$, we know that any $G$-invariant
probability measure on $M \atimes I$ that is equivalent to
$\mu \times \lambda$ must project to~$\mu$ on~$M$.)

\begin{notation}
 For each $m \in M$, let $\lambda_m$ be the Lebesgue
measure on the interval $m \times I$. Thus,
 $\lambda_m(E) = \lambda_m \bigl( E \cap (m \times I)
\bigr)$ for every Borel subset~$E$ of~$I$, and we have
 $\mu \times \lambda = \int_M \lambda_m \, d\mu(m)$.
 \end{notation}

 {\bf Assume} that  $(\mu \times \lambda) \bigl(
\phi^{-1}(0) \bigr) \neq 0$. (This will lead to a
contradiction.) Because $\mu$ is ergodic and
$\phi^{-1}(0)$ is $G$-invariant, we must have $\lambda_m
\bigl( \phi^{-1}(0) \bigr) \neq 0$, for a.e.~$m \in M$. By
discarding a set of measure~$0$, we may assume $\lambda_m
\bigl( \phi^{-1}(0) \bigr) \neq 0$ for
every $m \in M$. Then we may define $f \colon M \to I$ by
 $$ f(m) = \max\bigset{ t \in I}
 {\lambda_m \Bigl( \phi^{-1}(0) \cap \bigl(m \times
[0,t]\bigr) \Bigr)  = 0} .$$
 Replacing $M \atimes I$ with the invariant subset 
 $$ \bigl\{ \, (m,s) \in M \atimes I \mid f(m) \le
s \le 1 \, \bigr\} ,$$
 we may assume that 
 \begin{equation} \equationlabel{f=0}
 \mbox{$f$ is identically~$0$.}
 \end{equation}

\setcounter{step}{0}

\begin{step} \label{phinotvanish}
 There is some $\delta_0 > 0$, such that
 $ \mu \bigset{ m \in M}
 {\int_{m \times I'} |\phi|^2 \, d \lambda_m < \delta_0 }
 < \epsilon_2 $.
 \end{step}
 (To avoid confusion, we emphasize that the integral is over
$m \times I'$, not the entire interval $m \times I$.)
  Let $\chi$ be the characteristic function of $M \times
I'$. We will show, for every $g \in C$, that
 $\| \rho(g)\chi - \chi \| < \delta_2 \|\chi\|$
(see Claim~\ref{chi-alminv} below).
 From the definition of~$\delta_2$, this implies that there
is some $\psi \in \invariant$, such that $\|\psi\| =
\|\chi\|$ and $\| \chi - \psi \| \le \epsilon_2 \|\chi\|$.
Therefore,
 $$ \mu \{\, m \in M \mid \text{$\psi = 0$ a.e.\ on $m
\times I'$} \,\} \le \epsilon_2^2 < \epsilon_2 .$$
 From \pref{phi-1(0)}, we conclude that the same inequality
is true with $\phi$ in the place of~$\psi$. In other
words, we have
 $$ \mu \bigset{ m \in M}
 {\int_{m \times I'} |\phi|^2 \, d \lambda_m = 0 }
 < \epsilon_2 .$$
 Thus, the desired conclusion is obtained by taking
$\delta_0$ sufficiently small.

\begin{stepclaim} \label{chi-alminv}
 For each $g \in C$, we have
 $\| \rho(g)\chi - \chi \| < \delta_2 \|\chi\|$.
 \end{stepclaim}
 Let
 $$E = \{ m \in M \mid \forall s \in 2I', ~ 
 |\alpha(g^{-1},m)'(s) - 1| < \epsilon_3 \} .$$
 For $m \in E$, we have:
 $$ \Bigl| \bigl( \rho(g)\chi \bigr) (m,s) - \chi(m,s)
\Bigr| \le
 \begin{cases}
 \epsilon_3 & \text{if $(1+\epsilon_3)s \in I'$} \\
 2 &  \text{if $(1+\epsilon_3)s \notin I'$ and 
 $s/(1+\epsilon_3) \in I'$} \\
 0 &  \text{if $s/(1+\epsilon_3) \notin I'$} \\
 \end{cases}
 $$
 Therefore, 
 \begin{eqnarray}\label{int(E)chi2}
 \int_{E \times I} &&
 \bigl| \bigl( \rho(g)\chi \bigr) - \chi \bigr|^2
 \, d(\mu \times \lambda) \\
 &&\le \epsilon_3^2 \lambda(I')
 + 2^2 \left( (1+\epsilon_3) -
\frac{1}{1+\epsilon_3} \right) \lambda(I')
 + 0
 \le 9 \epsilon_3 \lambda(I') . \notag
 \end{eqnarray}
 
If $F$ is any subset of~$M$ with $\mu(F) < \epsilon_3$, then
$\int_{F \times I} \chi^2 \, d(\mu \times \lambda) <
\epsilon_3 \lambda(I')$, because $|\chi| \le 1$. Then we must also have 
$\int_{F \times I}\bigl( \rho(g)\chi \bigr)^2 \, d(\mu
\times \lambda) < \epsilon_3 \lambda(I')$, because $\rho$
is unitary, $g^{-1}(F \times I) = (g^{-1}F) \times I$, and
$\mu$ is $G$-invariant. In particular,
letting $F = M \setminus  E$ and using the triangle
inequality, we obtain
 \begin{equation} \equationlabel{int(M-E)chi2}
 \int_{(M \setminus E) \times I} \Bigl|
\bigl( \rho(g)\chi \bigr) - \chi \Bigr|^2
 \, d(\mu \times \lambda)
 < 4 \epsilon_3 \lambda(I') .
 \end{equation}

Combining \pref{int(E)chi2}
and~\pref{int(M-E)chi2} yields
 $$\| \rho(g)\chi - \chi \|^2
 \le 9 \epsilon_3 \lambda(I') + 
4 \epsilon_3 \lambda(I')
 = 13 \epsilon_3 \| \chi \|^2
 < \bigl( \delta_2 \|\chi\| \bigr)^2 .$$
 This completes the proof of the claim.

\begin{step}
 We obtain a contradiction.
 \end{step}
 Let $\chi'$ be the characteristic function of the
$G$-invariant set 
 $$X = \phi^{-1}(0) \cap
 \bigset{ (m,s) \in M \atimes I }
 {\int_0^s |\phi(m,t)|^2 \, d \lambda_m(t) < \delta_0 } .$$
 We have $\lambda_m(X) \neq 0$ for a.e.\ $m \in M$
\see{f=0}, so we may define a unit vector $\omega \in L^2(M
\atimes I)$ by
 $$ \omega(m,s) = \frac{\chi'(m,s)}{\lambda_m (X)^{1/2}} .$$
  We will show, for every $g \in C$, that
 $\| \rho(g)\omega - \omega \| < \delta_1 \|\omega\|$ (see
Claim~\ref{omega-alminv} below).
 Then the definition of~$\delta_1$ implies that $\omega$ is
not orthogonal to~$\invariant$. Thus, there is some $\psi
\in \invariant$, such that $\psi$ is not essentially~$0$
on~$X$. From~\pref{phi-1(0)}, we conclude that $\phi$ is not
essentially~$0$ on~$X$. Because $X \subset \phi^{-1}(0)$,
this is a contradiction.

\begin{stepclaim} \label{omega-alminv}
 For every $g \in C$, we have
 $\| \rho(g)\omega - \omega \| < \delta_1 \|\omega\|$.
 \end{stepclaim}
 Let 
 $$ E = \bigset{ m \in M }{ \forall s \in 2I', ~ 
 |\alpha(g^{-1},m)'(s) - 1| < \epsilon_3 }
 \cap  \bigset{ m \in M}
 {\int_{m \times I'} |\phi|^2 \, d \lambda_m \ge \delta_0 }
.$$
 By comparing the rightmost terms in the definitions of~$X$
and~$E$, we see that $X \cap (E \times I) \subset E \times
I'$. Thus, from the left term in the definition of~$E$, we
see that
 $|\alpha(g^{-1},m)'(s) - 1| < \epsilon_3$ for every $(m,s)
\in X \cap (E \times I)$. 
 Therefore, for $m \in E$, we have:
 $$ \Bigl| \bigl( \rho(g)\omega \bigr) (m,s) - \omega(m,s)
\Bigr| \le
 \begin{cases}
 {\epsilon_3} /
 {\lambda_m (X)^{1/2}}
 & \text{if $(m,s) \in X$} \\
 0 & \text{if $(m,s) \notin X$} \\
 \end{cases}
 $$
 Therefore, 
 \begin{equation} \equationlabel{int(E)omega}
 \int_{E \times I} 
 \Bigl| \bigl( \rho(g)\omega \bigr) - \omega \Bigr|^2
 \, d(\mu \times \lambda) \le \epsilon_3^2 < \epsilon_2
 . 
 \end{equation}
 
If $F$ is any subset of~$M$ with $\mu(F) < 2\epsilon_2$, then
$\int_{F \times I} \omega^2 \, d(\mu \times \lambda) <
2\epsilon_2$, because $\int \omega^2 \, d\lambda_m = 1$ for
every $m \in M$. Then we must also have
 $\int_{F \times
I}\bigl( \rho(g)\omega \bigr)^2 \, d(\mu \times \lambda) <
2\epsilon_2$, because $\rho$ is unitary, $g^{-1}(F \times I)
= (g^{-1}F) \times I$, and $\mu$ is $G$-invariant. In
particular, letting $F = M \setminus  E$ and using the
triangle inequality, we obtain
 \begin{equation} \equationlabel{int(M-E)omega}
 \int_{(M \setminus E) \times I} \Bigl| \bigl(
\rho(g)\omega \bigr) - \omega \Bigr|^2
 \, d(\mu \times \lambda)
 < 8 \epsilon_2 .
 \end{equation}

Combining \pref{int(E)omega} and~\pref{int(M-E)omega} yields
 $$\| \rho(g)\omega - \omega \|^2
 \le \epsilon_2 +  8 \epsilon_2 < \delta_1^2 
 = \bigl( \delta_1 \|\omega\| \bigr)^2 .$$
 This completes the proof of the claim.
 It also completes the proof of Theorem~\ref{Kazhdan-homeo}.
 \end{pf*}

\begin{rem}
 For a smooth manifold~$X$ and a point $x \in X$, let
$\Diff^1_{\Id}(X;x)$ be the group of
$C^1$~diffeomorphisms~$h$ of~$X$, such that $h(x) = x$ and
$Dh(x) = \Id$. It would be interesting to know whether
Theorem~\ref{Kazhdan-homeo} generalizes to the cocycles
$\alpha \colon G \times M \to \Diff^1_{\Id}(X;x)$, for
$\dim X > 1$.

It would also be interesting to know whether additional
smoothness on the cocycle~$\alpha$ yields additional
smoothness on the function that implements the cohomology
of~$\alpha$ to a trivial cocycle.
 \end{rem}

\section{Other versions of the main theorem}
\label{versions-section}

The assumption that $M$ is irreducible and ergodic is
stronger than is necessary in Theorems~\ref{GhysThm}
and~\ref{C1-isometry}. Namely, we may allow $G$ to be a
product of higher-rank normal subgroups whose ergodic
components are irreducible \see{irred-factors}. 
In particular, if no  simple factor of~$G$ has real rank
one, then there is no need for any ergodicity or
irreducibility assumption on~$M$ \see{no-rank1}.
 
There are also analogous results in the more general
situation where $G$ is allowed to be a product of semisimple
algebraic groups over local fields \see{Ghys-S-algebraic},
or a lattice in such a group (see~\ref{Ghys-Lattice}
and~\ref{Ghys-S-arithmetic}). Thus,
Theorem~\ref{GhysC1LatticeThm} generalizes to the
situation where $\Gamma$ is an $S$-arithmetic group
(see~\ref{Ghys-Slattice} and~\ref{Sarith-on-circle}).

The results also generalize to the case where $G$ has
$\PSL(2,\real)$ as a factor, but the conclusion must be
weakened (see~\ref{Ghys-SL2} and~\ref{Ghys-SL2cocycle}).

\begin{cor} \label{irred-factors}
 Let 
 \begin{itemize}
 \item $G$ be a connected, semisimple, real Lie group, such that
there is no continuous homomorphism from~$G$ onto
$\PSL(2,\real)$;
 \item $M$ be a Borel $G$-space with invariant probability
measure~$\mu$;
 \item $\alpha \colon G \times M \to \Homeo(\torus)$ be a
strict Borel cocycle; and
 \item $G_0,G_1,\ldots,G_r$ be connected,
closed, normal subgroups of~$G$, such that
 \begin{itemize}
 \item $G = G_0 G_1 \cdots G_r$,
 \item $G_0$ is compact, 
 \item $\Rrank(G_i) \ge 2$ for $i > 0$,
 \item $[G_i,G_j] = e$ for all $i$ and~$j$, and
 \item for each~$i > 0$, almost every ergodic component of
the action of~$G_i$ on~$M$ is irreducible.
 \end{itemize}
 \end{itemize}
 Then there is a $G$-invariant probability measure~$\nu$ on
$M \atimes \torus$, such that the projection of~$\nu$
to~$M$ is~$\mu$.

Furthermore, if $G$ has Kazhdan's property~$(T)$, and
$\alpha(G \times M) \subset \Diff^1(\torus)$, then $\nu$
can be taken to be equivalent to $\mu \times \lambda$.
 \end{cor}

\begin{pf}
 If $G = G_0$, then $G$ is compact, so we obtain a
$G$-invariant measure~$\nu$ on $M \atimes \torus$ simply by
averaging over~$G$.
 Thus, we may assume $r > 0$. Let
 \begin{itemize}
 \item  $G^* = G_0 G_1 \cdots G_{r-1}$;
 \item $\alpha^*$ be the restriction of~$\alpha$ to $G^*
\times M$;
 \item $\alpha_r$ be the restriction of~$\alpha$ to $G_r
\times M$; and
 \item $\mathcal{A} = 
 \bigl\{ \psi \colon M \to \Prob(\torus) \mid
 \text{$\psi$ is $\alpha^*$-equivariant} \,\}$.
 \end{itemize}
 By induction on~$r$, we may assume that there is a
$G^*$-invariant probability measure $\nu^*$ on $M
\mathbin{\times_{\alpha^*}} \torus$, such that the
projection of~$\nu^*$ to~$M$ is~$\mu$. Therefore,
$\mathcal{A}$ is nonempty.

Let $P$ be a minimal parabolic subgroup of~$G$. The space
of measurable functions from $G_r/(P \cap G_r)$
to~$\mathcal{A}$ is an affine $G_r$-space over~$G_r/(P \cap
G_r)$, so, because $P \cap G_r$ is amenable, there is a
$\beta$-equivariant Borel map $\Phi \colon G_r/(P \cap G_r)
\to \mathcal{A}$. Then, defining $\Psi \colon \bigl( G_r/(P
\cap G_r) \bigr) \times M \to \Prob(\torus)$ by $\Psi(x, m)
= \Phi(x)(m)$, we see that $\Psi$ is
$\alpha_r$-equivariant. 

For each $m \in M$, define $\Psi_m \colon G/P \to
\Prob(\torus)$ by $\Psi_m(x) = \Psi(x,m)$.
By assumption, each ergodic component of the action
of~$G_r$ on~$M$ is irreducible. Thus,
Theorem~\ref{Ghys-essconst} implies that $\Psi_m$ is
essentially constant, for a.e.~$m \in M$. Thus, $\Psi$
induces an essentially well-defined Borel map $\bar\Psi
\colon M \to \Prob(\torus)$. By construction, $\bar\Psi$ is
both $\alpha^*$-equivariant and $\alpha_r$-equivariant, so
$\bar\Psi$ is $\alpha$-equivariant. Therefore, $\nu =
\int_M \bigl( m \times \bar\Psi(m) \bigr) \, d\mu(m)$ is a
$G$-invariant measure on $M \atimes \torus$. By
construction, it projects to~$\mu$ on~$M$.
 \end{pf}

\begin{cor} \label{no-rank1}
 Let 
 \begin{itemize}
 \item $G$ be a connected, semisimple, real Lie group, such that
$G$ has no factors of real rank~one;
 \item $M$ be a Borel $G$-space with invariant probability
measure~$\mu$; and
 \item $\alpha \colon G \times M \to \Homeo(\torus)$ be a
strict Borel cocycle.
 \end{itemize}
 Then there is a $G$-invariant probability measure~$\nu$ on
$M \atimes \torus$, such that the projection of~$\nu$
to~$M$ is~$\mu$.

Furthermore, if $\alpha(G \times M) \subset
\Diff^1(\torus)$, then $\nu$ can be taken to be equivalent
to $\mu \times \lambda$.
 \end{cor}

\begin{cor} \label{Ghys-Lattice}
 Let 
 \begin{itemize}
 \item $G$ be a connected, semisimple, real Lie group, such that
 \begin{itemize}
 \item $\Rrank G \ge 2$, and 
 \item there is no continuous homomorphism from~$G$ onto
$\PSL(2,\real)$;
 \end{itemize}
 \item $\Gamma$ be an irreducible lattice in~$G$;
 \item $M$ be an irreducible ergodic $\Gamma$-space with
finite invariant measure~$\mu_M$; and
 \item $\alpha \colon \Gamma \times M \to \Homeo(\torus)$
be a Borel cocycle.
 \end{itemize}
 Then there is a $\Gamma$-invariant probability
measure~$\nu$ on $M \atimes \torus$, such that the
projection of~$\nu$ to~$M$ is~$\mu$.

Furthermore, if $\Gamma$ has Kazhdan's property~$(T)$, and
$\alpha(\Gamma \times M) \subset \Diff^1(\torus)$, then
$\nu$ can be taken to be equivalent to $\mu \times \lambda$.
 \end{cor}

\begin{pf*}{Proof {\rm (a standard argument)}}
 Let $\funddom \colon G/\Gamma \to G$ be a Borel section
(i.e., $\funddom(g \Gamma) \in g \Gamma$ for every $g
\Gamma \in G/\Gamma$), and define $\gamma \colon G \times
G/\Gamma \to \Gamma$ by
 $ \funddom(g x) = g \funddom(x) \gamma(g,x)^{-1}$
 for $g \in G$ and $x \in G/\Gamma$.
 Then $\gamma$ is a Borel cocycle for the action of~$G$
on~$G/\Gamma$.

Let $\hat M = \operatorname{Ind}_\Gamma^G(M) = (G/\Gamma)
\mathbin{\times_\gamma} M$ be the $G$-space induced
from the $\Gamma$-space~$M$. Then $\hat M$ is an
irreducible, ergodic $G$-space.
 Define a Borel cocycle $\hat\alpha \colon G \times \hat M
\to \Homeo(\torus)$ by
 $ \hat\alpha \bigl( g, (x,m) \bigr) = \alpha \bigl(
\gamma(g,x), m \bigr)$.

From Theorem~\ref{GhysThm}, we know that there is a
$G$-invariant probability measure~$\hat\nu$ on~$\hat M
\mathbin{\times_{\hat\alpha}} \torus$, such that $\hat\nu$
projects to~$\rho \times \mu$ on $(G/\Gamma) \times M$,
where $\rho$ is the $G$-invariant probability measure on
$G/\Gamma$.
 We may write
 $ \hat\nu = \int_{G/\Gamma} \hat\nu_x \, d \rho(x)$,
 where $\hat\nu_x \in \Prob(x \times M \times \torus)$.

 Fix a.e.\ $g \in G$. The measure $\hat\nu_{g\Gamma}$ is $g
\Gamma g^{-1}$-invariant, and projects to~$\mu$ on~$M$.
Define $\nu$ by $g_*(e\Gamma \times \nu) = \nu_{g\Gamma}$.
 \end{pf*}

\begin{rem}
 For $G$ as in Theorem~\ref{Ghys-S-algebraic}, the
definitions of irreducible lattice and irreducible action
given in~\S\ref{intro} (see~\ref{irredlattice}
and~\ref{irredaction}) must be modified to refer to
``non-discrete" normal subgroups instead of ``connected" normal
subgroups.
 \end{rem}

\begin{thm} \label{Ghys-S-algebraic}
 Let 
 \begin{itemize}
 \item $\mathcal{S}$ be a finite set of local fields
{\upshape(}not necessarily of characteristic
zero{\upshape)};
 \item $\algG_F$ be a connected, semisimple algebraic
group over~$F$, for each $F \in \mathcal{S}$;
 \item $G$ be a closed, cocompact, normal subgroup of
$\prod_{F \in \mathcal{S}} \algG_F(F)$;
 \item $M$ be an irreducible ergodic $G$-space with finite
invariant measure~$\mu$; and
 \item $\alpha \colon G \times M \to \Homeo(\torus)$ be a
strict Borel cocycle.
 \end{itemize}
 Assume
 \begin{enumerate}
 \renewcommand{\theenumi}{\alph{enumi}} 
 \item \label{Ghys-S-algebraic-Frank}
 $\sum_{F \in \mathcal{S}} \Frank{F}(\algG_F) \ge
2$; and
 \item \label{Ghys-S-algebraic-SL2}
 the identity component $G^\circ$ has no continuous
homomorphism onto $\PSL(2,\real)$.
 \end{enumerate}
  Then there is a $G$-invariant probability measure~$\nu$
on $M \atimes \torus$, such that the projection of~$\nu$
to~$M$ is~$\mu$.

Furthermore, if $G$ has Kazhdan's property~$(T)$, and
$\alpha(G \times M) \subset \Diff^1(\torus)$, then
$\nu$ can be taken to be equivalent to $\mu \times \lambda$.
 \end{thm}

\begin{pf*}{Remarks on the proof}
 Most of the arguments of Sections~\ref{ProveGhys}
and~\ref{Thurston-section} apply with only minor
changes. 

As a replacement for the Moore Ergodicity
Theorem~\ref{MooreErgodicityThm}, we note that the proof
of \cite[Thm.~7.2, p.~105]{MargulisBook} yields a version
of this result that applies to the general groups~$G$ under
consideration, in the special case where $H$ contains a
nontrivial split torus of~$G$. This suffices for our purposes.

When $F$ is nonarchimedean, we use Lemma~\ref{SL2inv(charp)} below
in place of Theorem~\ref{SL2inv}. (The proof of this lemma relies
on an argument of \'E.~Ghys \cite[pp.~219--220]{Ghys} that was not
needed in \S\ref{ProveGhys} or~\S\ref{Thurston-section}.) The
existence of a subgroup~$L$ satisfying the hypotheses of this
lemma follows from Lemma~\ref{littleSL2} below. (Because
$\SL(2,F)$ has no infinite, proper, normal subgroups
\cite[Cor.~2.3.2, p.~53]{MargulisBook}, we know that $\zeta \bigl(
\SL(2,F) \bigr) \subset G$.)
 \end{pf*}

\begin{lem}[{\cite[Prop.~3.1(13)]{TitsAbstractGrps}}]
\label{littleSL2}
 Let 
 \begin{itemize}
 \item $\algG$ be a semisimple algebraic group
over a field~$F$;
 \item $\algA$ be a maximal $F$-split torus
of~$\algG$; and
 \item $\alpha$ be an $F$-root of~$\algG$ {\upshape(}with
respect to~$\algA${\upshape)}, such that $2\alpha$ is not
an $F$-root.
 \end{itemize}
 Then there is a nontrivial $F$-homomorphism $\zeta \colon
\SL(2,\cdot) \to \algG$, such that
 $$ \zeta \begin{pmatrix}
 1 & x \\
 0 & 1
 \end{pmatrix}
 \in \algU_{\alpha} (F)
 , \qquad
 \zeta \begin{pmatrix}
 1 & 0 \\
 x & 1
 \end{pmatrix}
 \in \algU_{-\alpha} (F)
 , \qquad 
 \text{and}
 \qquad
 \zeta \begin{pmatrix}
 a & 0 \\
 0 & a^{-1}
 \end{pmatrix}
 \in \algA (F)
 , $$
 for all $x \in F$ and $a \in F \setminus \{0\}$.
 \end{lem}

We will use the following elementary observation in the
proof of Lemma~\ref{SL2inv(charp)}.

\begin{lem} \label{SL2(nonreal)}
 Let
 \begin{itemize}
 \item $F$ be a local field;
 \item $L = \SL(2,F)$; 
 \item $P$ be a proper parabolic subgroup of~$L$; and
 \item $a,b,c$ be three distinct elements of~$L/P$.
 \end{itemize}
 If $F \neq \real$, then there exist
 $y_0,\ldots,y_n \in T$, such that
 $y_0 = b$, $y_n = a$, and $(y_{i-1},y_i,c) \in
L(a,b,c)$ for $i = 1,\ldots,n$.
 \end{lem}

\begin{pf}
 It is well known that there is an identification of $L/P$ with $F
\cup \{\infty\}$, so that we have the standard action of~$L$ on $F
\cup \{\infty\}$ by linear-fractional transformations. Then,
because $\GL(2,F)$ is triply transitive on $F
\cup \{\infty\}$ and normalizes~$L$,
we may assume $a = 0$, $b = 1$, and $c = \infty$. Because $F \neq
\real$, we may choose $t_0,\ldots,t_n \in F$ such that $t_0 = 1$
and $t_0^2 + \cdots + t_n^2 = 0$. Let $y_i = t_0^2 + \cdots +
t_i^2$. Then
 $$ (y_{i-1},y_i,c) = 
 \begin{pmatrix}
 t_{i-1} & y_{i-1}/t_{i-1} \\
 0   & 1/t_{i-1}
 \end{pmatrix}
 ( 0 , 1 , \infty)
 \in L(a,b,c)
 ,$$
 as desired.
 \end{pf}

\begin{lem} \label{SL2inv(charp)}
 Let $G$, $M$, and~$\alpha$ be as in
Theorem~\ref{Ghys-S-algebraic}, and let $P = G \cap
\prod_{F \in \mathcal{S}} \mathord{\mathbf{P}}_F(F)$, where
$\mathord{\mathbf{P}}_F$ is a minimal parabolic subgroup
of~$\algG_F$, for each $F \in \mathcal{S}$.

 Suppose $\Psi \colon (G/P) \times M \to \torus_k$ is an
$\alpha$-equivariant Borel map. For each $m \in M$, define
$\Psi'_m \colon G \to \torus_k$ by $\Psi'_m(g) =
\Psi(gP,m)$.

 Suppose $F \neq \real$, and $L$ is a closed subgroup
of~$G$, such that
 \begin{itemize}
 \item \label{SL2inv(charp)-noncpct}
 $C_P(L)$ contains a nontrivial split torus of~$G$;
 \item \label{SL2inv(charp)-SL2}
 $L$ is the image of an $F$-morphism $\zeta \colon
\SL(2,\cdot) \to \algG_F$ with finite kernel; and
 \item \label{SL2inv(charp)-P}
 $\zeta^{-1}(L \cap P)$ is a parabolic subgroup of~$\SL(2,F)$.
 \end{itemize}
 Then $\Psi'_m$ is essentially right $L$-invariant, for
a.e.~$m \in M$.
 \end{lem}

\begin{pf}
 Define $X$ as in the proof of Theorem~\ref{SL2inv}, and fix some
$(a,b,c) \in X$. It follows from Lemma~\ref{SL2(nonreal)} that
there exist $y_0,\ldots,y_n \in G/P$, such that
 $y_0 = b$, $y_n = a$, and $(y_{i-1},y_i,c) \in G (a,b,c)$ for $i
= 1,\ldots,n$.

\setcounter{case}{0}

\begin{case} \label{SL2inv(charp)pf-k=1}
 Assume that $k = 1$, and that $\Psi'_m(x_1)$, $\Psi'_m(x_2)$, and
$\Psi'_m(x_3)$ are distinct, for almost all $m \in M$ and
$(x_1,x_2,x_3) \in G(a,b,c)$.
 \end{case}
 From the Moore Ergodicity Theorem, we know that $G$ is ergodic on
$G(a,b,c) \times M$. Thus, for almost every $m \in M$ and $g \in
G$, we conclude that 
 $$ \bigl( \Psi'_m( g y_{i-1} ), \Psi'_m( g y_i ), \Psi'_m(
g c ) \bigr)
 \text{\ \ and\ \ }
 \bigl( \Psi'_m( g a ), \Psi'_m( g b ), \Psi'_m( g c )
\bigr) $$
  have the same orientation, for $i = 1,\ldots,n$.
 Thus, by induction, we see that
 $$ \bigl( \Psi'_m( g y_0 ), \Psi'_m( g y_i ), \Psi'_m( g c )
\bigr)
 \text{\ \ and\ \ }
 \bigl( \Psi'_m( g a ), \Psi'_m( g b ), \Psi'_m( g c )
\bigr) $$
 have the same orientation.
 In particular, by letting $i = n$, we see that
 $$ \bigl( \Psi'_m( g b ), \Psi'_m( g a ), \Psi'_m( g c )
\bigr)
 \text{\ \ and\ \ }
 \bigl( \Psi'_m( g a ), \Psi'_m( g b ), \Psi'_m( g c )
\bigr) $$
 have the same orientation.
 This is a contradiction.

\begin{case}
 The general case.
 \end{case}
 If $\Psi_m'$ is not essentially right $L$-invariant, then the
argument of \cite[pp.~219--220]{Ghys} shows that, after
replacing~$\Psi$ with a different $\alpha$-equivariant Borel
function, we may assume, for almost every $m \in M$ and
$(x_1,x_2,x_3) \in G(a,b,c)$, that $\Psi_m'(x_1)$ is disjoint from
$\Psi_m'(x_2)$, and the sets $\Psi_m'(x_1)$ and $\Psi_m'(x_2)$
alternate around the circle. 

Now, if three
pairwise disjoint $k$-element subsets $B_1,B_2,B_3$
of~$\torus$ are pairwise alternating around the circle, we
say that $(B_1,B_2,B_3)$ is 
\emph{positively oriented} if there is a positively
oriented arc of the circle from a point of~$B_1$ to a
point of~$B_2$ that does not contain any point of~$B_3$
\cite[bot.\ of p.~220]{Ghys}. This relation has the
properties of a circular order, so we obtain a contradiction by
applying the same argument as in Case~\ref{SL2inv(charp)pf-k=1}.
 \end{pf}

The proof of Corollary~\ref{Ghys-Lattice} yields the
following as a corollary of Theorem~\ref{Ghys-S-algebraic}.

\begin{cor} \label{Ghys-S-arithmetic}
 Let 
 \begin{itemize}
 \item $G$ be as in Theorem~\ref{Ghys-S-algebraic}
{\upshape(}including assumptions
\pref{Ghys-S-algebraic-Frank}
and~\pref{Ghys-S-algebraic-SL2}{\upshape)}; 
 \item $\Gamma$ be an irreducible lattice in~$G$;
 \item $M$ be an irreducible ergodic $\Gamma$-space
with finite invariant measure~$\mu$; and
 \item $\alpha \colon \Gamma \times M \to \Homeo(\torus)$ be
a strict Borel cocycle.
 \end{itemize}
 Then there is a $\Gamma$-invariant probability
measure~$\nu$ on $M \atimes \torus$, such that the
projection of~$\nu$ to~$M$ is~$\mu$.

Furthermore, if $\Gamma$ has Kazhdan's property~$(T)$, and
$\alpha(\Gamma \times M) \subset \Diff^1(\torus)$, then
$\nu$ can be taken to be equivalent to $\mu \times \lambda$.
 \end{cor}

The following generalization of Theorems~\ref{GhysC1LatticeThm}
and~\ref{GhysHomeoLattThm} is the special case of
Corollary~\ref{Ghys-S-arithmetic} in which $M$ is a single point.
This result is essentially due to M.~Burger and N.~Monod
\cite{BurgerMonod-1, BurgerMonod-2, BurgerMonod-3}, but a few
isolated cases are not covered by their theorems. (On the other
hand, some of their results apply in a more general setting where
$G$ need not be a product of algebraic groups, or Lie groups.) In
the final conclusion of this corollary, we do not assume $\Gamma$
has Kazhdan's property~$(T)$, because
Thurston's Theorem~\ref{ReebThurston} can be applied if $\Gamma$ is
finitely generated. Furthermore, this restriction to finitely
generated lattices may be superfluous: we do not know an example
of an irreducible lattice in such a group that is not finitely
generated.

\begin{cor} \label{Ghys-Slattice}
 Let 
 \begin{itemize}
 \item $G$ be as in Theorem~\ref{Ghys-S-algebraic}
{\upshape(}including assumptions
\pref{Ghys-S-algebraic-Frank}
and~\pref{Ghys-S-algebraic-SL2}{\upshape)}; and
 \item $\Gamma$ be an irreducible lattice in~$G$.
 \end{itemize}
 Then every continuous action of~$\Gamma$ on~$\torus$ has
a finite orbit.

Furthermore, if $\Gamma$ is finitely generated, then every
homomorphism from~$\Gamma$ to $\Diff^1(\torus)$ has finite
image.
 \end{cor}

Ghys' Theorems~\ref{GhysC1LatticeThm}
and~\ref{GhysHomeoLattThm} are essentially the special
case of the following corollary in which $E$ is a number
field and $S$ consists only of the infinite places. (More
generally, if $E$ is a number field and
 $\sum_{s \in S_\infty}\Frank{E_s}(\algG) \ge 2$,
 then Conclusion~\pref{Sarith-on-circle-diffeo} of this
corollary is a consequence of Ghys' Theorem. Namely,
Theorem~\ref{GhysC1LatticeThm} applies to the subgroup
$\algG(\mathcal{O})$ of~$\Gamma$, and then the Margulis
Finiteness Theorem \cite[Thm.~IV.4.10,
p.~167]{MargulisBook} implies that the image of~$\Gamma$
is finite.) Note that Assumption~\pref{Sarith-on-circle-localrank}
implies $\Gamma$ is finitely generated \cite[Thm.~III.5.7(c),
p.~131]{MargulisBook}, \cite[Thm.~1a]{Behr}.

\begin{cor} \label{Sarith-on-circle}
 Let 
 \begin{itemize}
 \item $E$ be a global field;
 \item $S$ be a nonempty, finite set of places of~$E$,
including all of the infinite places;
 \item $\algG$ be a connected, almost simple algebraic
group over~$E$;
 \item $\mathcal{O}(S)$ be the ring of $S$-integers in~$E$;
and
 \item $\Gamma$ be a finite-index subgroup of~$\algG \bigl(
\mathcal{O}(S) \bigr)$.
 \end{itemize}
 Assume
 \begin{enumerate}
 \renewcommand{\theenumi}{\alph{enumi}} 
 \item \label{Sarith-on-circle-localrank}
 $\sum_{s \in S} \Frank{E_s}(\algG) \ge 2$; and
 \item for each archimedean $s \in S$, there is no
continuous homomorphism from $\algG(E_s)^\circ$ onto
$\PSL(2,\real)$.
 \end{enumerate}
 Then 
 \begin{enumerate}
 \renewcommand{\theenumi}{\alph{enumi}} 
 \item every continuous action of~$\Gamma$ on~$\torus$ has
a finite orbit; and 
 \item \label{Sarith-on-circle-diffeo}
 every homomorphism from~$\Gamma$ to $\Diff^1(\torus)$ has
finite image.
 \end{enumerate}
 \end{cor}

The following theorem is the main result of
\cite[\S7]{Ghys}, although it was not stated explicitly.

\begin{thm}[{Ghys \cite[\S7]{Ghys}}] \label{Ghys-SL2xSL2}
 Let
 \begin{itemize}
 \item $G$ be a connected Lie group that is locally isomorphic to
$\SL(2,\real)^n$, for some $n > 0$;
 \item $\Gamma$ be a countable group;
 \item $\phi \colon \Gamma \to \Homeo_+(\torus)$ and $\iota
\colon \Gamma \to G$ be homomorphisms;
 \item $P$ be a parabolic subgroup of~$G$; and
 \item $\Psi \colon G/P \to \torus_k$ be a
$\Gamma$-equivariant Borel map, for some $k \ge 1$.
 \end{itemize}
 If $\iota(\Gamma)$ is ergodic on $G/H$, for every closed,
noncompact subgroup~$H$ of~$G$, then either
$\phi(\Gamma)$ has a finite orbit, or there is a
semiconjugacy as described in
Corollary~\fullref{Ghys-SL2}{semiconj} below.
 \end{thm}

Although Theorem~\ref{Ghys-SL2xSL2} assumes that $G$ is connected,
an examination of the proof shows that if $G$ is a real algebraic
group, then it holds under the weaker assumption that $G$ is
\emph{Zariski} connected. This yields the following generalization
of Corollary~\ref{Ghys-Slattice} that allows $\PSL(2,\real)$ as a
factor of~$G$. This generalization was proved by \'E.~Ghys
\cite[Thm.~1.2]{Ghys} for $\mathcal{S} \subset \{\real,
\complex\}$. To justify the stronger conclusion when $\phi(\Gamma)
\subset \Diff^2(\torus)$, see \cite[Prop.~10.2]{Ghys}.

\begin{cor} \label{Ghys-SL2}
 Let 
 \begin{itemize}
 \item $G$ be as in Theorem~\ref{Ghys-S-algebraic}, except
that we do not assume~\pref{Ghys-S-algebraic-SL2}
{\upshape(}although we do
assume~\pref{Ghys-S-algebraic-Frank}{\upshape)};
 \item $\Gamma$ be an irreducible lattice in~$G$; and
 \item $\phi \colon \Gamma \to \Homeo_+(\torus)$ be a homomorphism.
 \end{itemize}
 Then either
 \begin{enumerate}
 \item $\phi(\Gamma)$ has a finite orbit; or
 \item \label{Ghys-SL2-semiconj}
 the restriction of~$\phi$ to~$\Gamma$ is semiconjugate to
a finite cover of the composition of the following:
 \begin{enumerate}
 \item the inclusion of~$\Gamma$ into~$G$; 
 \item a continuous surjection $G \to \PSL(2,\real)$; and
 \item the standard action of $\PSL(2,\real)$ on~$\torus$ by
linear-fractional transformations.
 \end{enumerate}
 \end{enumerate}
 Furthermore, if $\phi(\Gamma) \subset \Diff^2(\torus)$,
then any semiconjugacy as in~\pref{Ghys-SL2-semiconj} above is
actually a topological conjugacy.
 \end{cor}

For completeness, we state the following generalization of
Theorem~\ref{Ghys-S-algebraic}. Its proof is completed by
translating \cite[\S7]{Ghys} in a straightforward way from the
setting of homomorphisms of lattices to the setting of Borel
cocycles for ergodic $G$-actions. 

\begin{cor} \label{Ghys-SL2cocycle}
 Let 
 $G$, $M$, $\mu$, and~$\alpha$ be as in
Theorem~\ref{Ghys-S-algebraic}, except that we do not
assume~\pref{Ghys-S-algebraic-SL2} {\upshape(}although we
do assume~\pref{Ghys-S-algebraic-Frank}{\upshape)}.
 Assume $\alpha(g,m)$ is orientation preserving, for all $g \in G$
and $m \in M$ \cf{GhysPf-orpres}.

 Then there is a probability
measure~$\nu$ on $M \atimes \torus$, such that the
projection of~$\nu$ to~$M$ is~$\mu$, and either
 \begin{enumerate}
 \item $\nu$ is $G$-invariant; or 
 \item there exist
 \begin{enumerate}
 \item a continuous surjection $\tau \colon G \to
\PSL(2,\real)$; and
 \item a $G$-equivariant, measure-preserving function 
  $$f \colon (M \atimes \torus, \nu) \to (M \ttimes \torus, \mu
\times \text{\upshape Leb}) ;$$
 \end{enumerate}
 such that $f$ is of the form $f(m,t) = \bigl( m, f_m(t) \bigr)$,
where, for a.e.~$m \in M$, 
 \begin{itemize}
 \item $f_m \colon \torus \to \torus$ is continuous, and 
 \item any continuous lift $\tilde f_m \colon \real \to \real$ is
increasing.
 \end{itemize}
 \end{enumerate}
 \end{cor}

\end{article}

\end{document}